\UseRawInputEncoding
\documentclass[12pt,a4paper]{article}

\usepackage{ucs}
\usepackage[usenames,dvipsnames]{xcolor}
\usepackage{graphics}
\usepackage{tkz-tab}
\usepackage{caption}
\usepackage{latexsym}
\usepackage{amssymb}
\usepackage{amsmath}
\usepackage{subcaption}
\usepackage{amsthm}
\newtheorem{theorem}{Theorem}[section]
\newtheorem{corollary}[theorem]{Corollary}
\newtheorem{lemma}[theorem]{Lemma}
\newtheorem{proposition}{Proposition}[section]
\newtheorem{definition}{Definition}[section]
\newtheorem{remark}{Remark}[section]
\usepackage{amsfonts,amssymb,eucal,amsmath}
\usepackage{bbold}

\title{On the finiteness  of    four-body central configurations }
\author{{ Xiang Yu $^1$\footnote{Email:xiang.zhiy@foxmail.com,  xiang.zhiy@gmail.com}}, {Shuqiang Zhu$^2$\footnote{corresponding author, ORCID: 0000-0002-6216-2326, Email:zhusq@swufe.edu.cn}} \\
\small \it School of Mathematics, Southwestern
University of Finance and Economics, \\
\small \it Chengdu 611130, China}

\date{}

\begin{document}
\maketitle

\begin{abstract}
The number of  central configurations in the four body problem  was proved to be finite, first by Hampton and Moeckel, then by Albouy and Kaloshin, when the masses are all positive. We prove that the four-body central configurations are finite for any four nonzero masses.  
\end{abstract}

\textbf{2020AMS Subject Classification}: {  70F10 \and 70F15 \and 	37Nxx }.

 ~~~~~~~~~~~~~~~~~~~~~~~~~~~~~~\textbf{Statements and Declarations}

 The authors declare that they have no conflicts of interest.

 Data sharing not applicable to this article as no datasets were generated or analysed during the current study.

 \tableofcontents

\section{Introduction}

We consider the finiteness of  central configurations in the four-body problem. Our approach  is to extend the method due to Albouy and Kaloshin. They first used this method to show that the number of central configurations in the four-body problem is finite if the masses are all positive, and that the number of  central configurations in the planar five-body problem is finite, perhaps except if  the 5-tuple of  positive masses belongs to a given codimension 2 subvariety of the mass space. We focus on the four-body case and obtain the finiteness result without regard to the sign of the masses, which settles the open question raised by Hampton and Moeckel in \cite{Hampton2006Finiteness}.

The $N$-body problem consists of describing the  complete behavior of solutions of the Newton's equations,  $m_k \ddot{{\mathbf{r}}}_k =\frac{\partial U}{\partial \mathbf {r}_k},   k=1, 2, \cdots, N$. Started by Newton and considered by many great mathematicians through the centuries, the problem for $N>2$ remains largely unsolved.  Central configurations arose in the study of the $N$-body problem. By definition, they are just special arrangement of the particles.  However, they play an important role in the study of the dynamics of the $N$-body problem. For instance, they are related to the homographic solutions, the analysis of collision orbits, and the bifurcation of integral manifold (cf. \cite{ Sma70-2, wintner1941analytical}).

The so-called  hypothesis on finiteness of central configurations  was
proposed by Chazy \cite{Chazy} and Wintner \cite{wintner1941analytical}, and was listed by Smale as the sixth problem  on
his list of problems for the 21-st  century \cite{Sma98}: Is the number of relative equilibria (planar central configurations) finite, in the $N$-body problem of celestial mechanics, for any choice of positive real numbers $m_1,\cdots,m_N$
as the masses?

Euler and Lagrange have solved this finiteness question when  $N=3$. In 2005, Hampton and Moeckel \cite{Hampton2006Finiteness} studied the case of $N=4$ , and showed that for any 4-tuple of positive masses the upper bound of numbers of central configurations  is 8472.  Albouy and  Kaloshin \cite{Albouy2012Finiteness} answered positively the question in 2012 for $N=5$,   except perhaps if the 5-tuple of positive masses satisfies two given polynomial conditions. We refer the reader to the excellent review on this problem by Hampton, Moeckel \cite{Hampton2006Finiteness} and Albouy, Kaloshin \cite{Albouy2012Finiteness}.

The outline of the method of  Albouy and Kaloshin is  the following: The central configurations are solutions of the  algebraic system  \eqref{equ:sys-0} in the real domain.   By introducing an efficient set of coordinates, Albouy and Kaloshin embedded system \eqref{equ:sys-0} into the system \eqref{equ:sys-3}. The solutions of  
 \eqref{equ:sys-3} in the complex domain are called \emph{ normalized central configurations}. If there are infinitely many normalized central configurations, Albouy and Kaloshin could extract a sequence of them approaching some singularities, which is called a \emph{singular sequence}. 
 
 These singular sequences are classified according to their limiting behaviors.  This classification was done by two-colored diagrams, which consists of colored strokes and circles and each circle represents a mass. For $N=4$, they found that the number of diagrams is 5, while the number increases to 16 for $N=5$.  
 
 They noted that for each diagram, or equivalently each class of singular sequences, there is some algebraic constraint on the masses, and another important fact that the existence of one class of singular sequence generally implies the existence of other classes. 
 Hence, the existence of infinitely many normalized central configurations would imply several algebraic constraints on the masses. The finiteness is asserted by showing that these constraints have no solutions.   In particular,  this method allowed to reduce  significant amount of computations.

 We study the finiteness question of   central configurations for any 4-tuple of nonzero masses.  It is worth the effort based on two reasons. Firstly, 
  negative masses shed light on  the role of masses in the N-body problem.  
 If  negative masses are allowed in the planar five-body problem, then Roberts \cite{roberts1999}  has constructed  a continuum of central configurations.   It is natural to ask   whether such  a continuum   is possible for the four-body problem if negative masses are allowed,  and  this is an open question  raised by Hampton and Moeckel in \cite{Hampton2006Finiteness}. Secondly, the  method  introduced by Albouy and Kaloshin   is  powerful, and the general idea could be implemented for  investigating other finiteness problem where  the restriction on parameters is weaker. Our study shows its potential in this aspect.   
 

In fact, negative masses, or even complex masses are not excluded in  the study of Albouy and Kaloshin   in \cite{Albouy2012Finiteness} although their main results are made for the positive masses. However, the notion of centers of mass of clusters is employed  in their study, and so it is presumed  that   \emph{the subset of bodies has total mass nonzero}. This hypothesis must be discarded in our study since we are considering the problem for all nonzero masses. 
  Thus, we shall proceed somewhat differently. In particular, we need to modify or abandon some rules in the classification of the singular sequences, see Section \ref{sec:rule}. As a result,  more diagrams appear. The number of diagrams  is now 14. Generally, the corresponding constraints on masses are more subtle,  even for the classes that also appear in Albouy and Kaloshin's  classification.  We carry out the proof separately for the case of the nonzero total mass and the case of zero total mass.

The results proved here are stated in the following
\begin{theorem}\label{thm:main}
Suppose that $m_1, m_2, m_3, m_4$ are real  and nonzero,  then   there are  finitely many  central configurations  in the   four-body
problem.
\end{theorem}

We remark that
\begin{theorem}\label{thm:main-m0}
 	Suppose that $m_1, m_2, m_3, m_4$ are real   and nonzero. 
 	Let  $\alpha$ be  the unique positive root of  the algebraic equation (\ref{xy}).   A complex continuum of relative equilibria  exists only for the two groups of masses: $m_1= m_2=- m_3=- m_4$,  $m_1= m_2,  m_3=-\alpha^2 m_1,m_4=-(2-\alpha^2) m_1$.  Numerically,   $\alpha\approx 1.2407$. 
 \end{theorem}

The paper is structured as follows. In Section \ref{sec:pre}, we introduce  some notations, definitions, the coordinate system and the notation of  singular sequences due to Albouy and  Kaloshin.  In Section \ref{sec:rule}, we  modify the tools  to classify the diagram corresponding to  singular sequences.  In Section \ref{sec:list},  we list  all possible diagrams for the four-body problem.  In Section \ref{sec:dia},  we  obtain the constraints on the masses corresponding to each diagram.  In Section \ref{sec:pr1} and \ref{sec:pr3}, we  outline  the proofs of  Theorem \ref{thm:main} and Theorem \ref{thm:main-m0}. Some heavy  computations in the proofs are postponed  in the Appendix.

 \section{Preliminaries}\label{sec:pre}

In this section, we recall,  among others notations and definitions,  the coordinate system and the notation of  singular sequences due to Albouy and  Kaloshin  \cite{Albouy2012Finiteness}.

  A central configuration is a solution of  the following  algebraic equations 
\begin{equation}\label{equ:sys-00}
\lambda  \mathbf{r}_k =\sum_{1 \leq j \leq N, j \neq k} \frac{m_j(\mathbf{r}_j-\mathbf{r}_k)}{|\mathbf{r}_j-\mathbf{r}_k|^3}, ~~~~~~~~~~~~~~~k=1, 2, \cdots, N.
\end{equation} 
Since we will consider the four-body problem which includes negative masses, 
the multiplier $\lambda$ can be positive, negative and even zero.   We will not study the case corresponding to $\lambda =0$. 
 Equations \eqref{equ:sys-00} is invariant under rotations and dilations in the plane. 
Hence,  we may assume that the multiplier equals $1$ or $-1$. 
Set $\mathbf{r}_k=(x_k, y_k) \in \mathbb{R}^2$, $k=1, ..., N$.  We rewrite equations \eqref{equ:sys-00} as 

  \begin{align}
 \label{equ:sys-0}
\sigma \begin{pmatrix}
 x_1\\y_1
 \end{pmatrix}&= m_2 r_{12}^{-3}\begin{pmatrix}
 x_{21}\\y_{21}
 \end{pmatrix}+m_3 r_{13}^{-3}\begin{pmatrix}
 x_{31}\\y_{31}
 \end{pmatrix}+...+m_N r_{1N}^{-3}\begin{pmatrix}
 x_{N1}\\y_{N1}
 \end{pmatrix}\\
 \notag  \sigma \begin{pmatrix}
 x_2\\y_2
 \end{pmatrix}&= m_1 r_{12}^{-3}\begin{pmatrix}
 x_{12}\\y_{12}
 \end{pmatrix}+m_3 r_{23}^{-3}\begin{pmatrix}
 x_{32}\\y_{32}
 \end{pmatrix}+...+m_N r_{2N}^{-3}\begin{pmatrix}
 x_{N2}\\y_{N2}
 \end{pmatrix}\\
 \notag & ...\\
 \notag  \sigma \begin{pmatrix}
 x_N\\y_N
 \end{pmatrix}&= m_1 r_{1N}^{-3}\begin{pmatrix}
 x_{1N}\\y_{1N}
 \end{pmatrix}+...+m_{N-1} r_{(N-1) N}^{-3}\begin{pmatrix}
 x_{(N-1)N}\\y_{(N-1)N}
 \end{pmatrix},
 \end{align}
 where $\sigma=\pm1, x_{kl}= x_l-x_k$, $y_{kl}=y_l-y_k$ and $r_{kl}=(x_{kl}^2+y_{kl}^2)^{1/2}>0$. Denote  by $f_k\in \mathbb{R}^2$, $k=1, ..., N$, the right-hand sides of the equations. System \eqref{equ:sys-0} can be abbreviated as $\sigma \mathbf{r}_k =f_k,  k=1, ..., N. $

We first embed the system (\ref{equ:sys-0}) above into a polynomial system in $\mathbb{C}^{2N}\times {\mathbb{C}^{N(N-1)/2}}$:
  \begin{align}
\label{equ:sys-1}
\begin{pmatrix}
x_1\\y_1
\end{pmatrix}&= m_2 \delta_{12}^{3}\begin{pmatrix}
x_{21}\\y_{21}
\end{pmatrix}+m_3 \delta_{13}^{3}\begin{pmatrix}
x_{31}\\y_{31}
\end{pmatrix}+... \\
\notag  \begin{pmatrix}
x_2\\y_2
\end{pmatrix}&= m_1 \delta_{12}^{3}\begin{pmatrix}
x_{12}\\y_{12}
\end{pmatrix}+m_3 \delta_{23}^{3}\begin{pmatrix}
x_{32}\\y_{32}
\end{pmatrix}+...\\
\notag & ...\\
\notag  &\delta_{12}^2(x_{12}^2+y_{12}^2)=1\\
\notag  &\delta_{13}^2(x_{13}^2+y_{13}^2)=1\\
\notag & ...\\
\notag &y_{12}=0.
\end{align}

 \begin{definition}[Normalized central configuration]\label{normalizedcentralconfiguration}
 	A normalized central configuration is a solution   of (\ref{equ:sys-1}). A real
 	normalized central configuration is a normalized central configuration such that
 	$(x_k,y_k)\in \mathbb{R}^2$ for any $k=1,2,\cdots,N$. A positive normalized central configuration
 	is a real normalized central configuration such that $\delta_{jk}=\pm 1/\sqrt{x_{jk}^2+y_{jk}^2}$
 	is
 	positive for any $j, k,j\neq k$
 \end{definition}

\begin{definition}The following quantities are defined:
\begin{center}
$\begin{array}{cc}
  \text{Total mass} & m =\sum_{j=1}^{N}m_j  \\
  \text{Potential function} & U =\sum_{1\leq j<k\leq N}\frac{m_jm_k}{r_{jk}}=\sum_{1\leq j<k\leq N}{m_jm_k}{\delta_{jk}}  \\
 \text{ Moment of mass }& M =\sum_{j=1}^{N}m_j \mathbf{r}_j \\
 \text{ Moment of inertia }& I =\sum_{j=1}^{N}m_j |\mathbf{r}_j|^2=\sum_{j=1}^{N}m_j (x_j^2+y_j^2)^2
\end{array}$
\end{center}

\end{definition}
Then it is easy to see that
\begin{equation}\label{mIr}
m I =\sum_{1\leq j<k\leq N}m_jm_k r_{jk}^2\triangleq S.
\end{equation}

\begin{lemma}[\cite{Albouy2012Finiteness}]\label{lem:domin}
	Let $X$ be a closed algebraic subset of $\mathbb{C}^\mathcal{N}$ and $f: \mathbb{C}^\mathcal{N} \mapsto \mathbb{C}$ be a polynomial. Either the image $f(X)\in \mathbb{C}$ is a finite set, or it it the complement of a finite set. In the second case one says that $f$ is dominating.
\end{lemma}

\begin{lemma}[\cite{Albouy2012Finiteness}]\label{lem:finitepotential}
	Consider the closed algebraic subset $\mathcal A\in \mathbb{C}^{2N} \times \mathbb{C}^{N(N-1)/2}$ defined by system \eqref{equ:sys-1} and the polynomial functions $U$, $I$ and $S$ on it. Then $U(\mathcal A)$, $I(\mathcal A)$ and $S(\mathcal A)$ are three finite sets.
\end{lemma}

 \subsection{Complex central configurations}

For notation convenience,  we may use  the variables $r_{jk}=1/\delta_{jk}$  instead of the $\delta_{jk}$'s.  But  keep in mind that the variables of system \eqref{equ:sys-1} are $x_k,y_k$ and $\delta_{jk}$.

Set $z_k = x_k + \mathbf{i}y_k$ and
$w_k = x_k - \mathbf{i}y_k$. If $x_k$'s and $y_k$'s are real,   the $z_k$'s form this configuration
in the complex plane, and  the $w_k$'s is the reflection 
of the configuration.  We have $x_k^2+y_k^2=z_kw_k$ 
and $x_{jk}^2+y_{jk}^2=z_{jk}w_{jk}$.  System \eqref{equ:sys-1} becomes

   \begin{align}
 \label{equ:sys-2}
z_1&= m_2 z_{21}^{-1/2} w_{21}^{-3/2} +m_3 z_{31}^{-1/2} w_{31}^{-3/2} +... \\
\notag  w_1&= m_2 z_{21}^{-3/2} w_{21}^{-1/2} +m_3 z_{31}^{-3/2} w_{31}^{-1/2} +... \\
\notag &...\\
\notag z_{12}&=w_{12}.
 \end{align}

 \indent\par

 	\emph{The word ``distance" will mean  $r_{jk}=\sqrt{z_{jk}{w_{jk}}}$. 	The $z_{jk}$'s (respectively the $w_{jk}$'s) in the complex plane will be called 
 	 $z$-separation (respectively $w$-separation). }

Set $Z_{jk}=z_{jk}^{-1/2} w_{jk}^{-3/2}$ and $W_{jk}=z_{jk}^{-3/2} w_{jk}^{-1/2}$, then system \eqref{equ:sys-2} becomes
\begin{align}
 \label{equ:sys-3}
z_1&= m_2 Z_{21} +m_3 Z_{31} +...+m_N Z_{N1}, \\
\notag  w_1&= m_2 W_{21} +m_3 W_{31} +...+m_N W_{N1}, \\
\notag &...\\
\notag z_N&= m_1 Z_{1N} +m_2 Z_{2N} +...+m_{N-1} Z_{(N-1)N}, \\
\notag  
w_N&= m_1 W_{1N} +m_2 W_{2N} +...+m_{N-1} W_{(N-1)N},\\
\notag z_{12}&=w_{12}.
 \end{align}
Here 
\begin{center}
$r_{jk}=r_{kj}=1/\sqrt[4]{Z_{jk}W_{jk}}$, 

$Z_{jk}=z_{jk}/r_{jk}^3$, $W_{jk}=w_{jk}/r_{jk}^3$, $Z_{jk}=-Z_{kj}$, $W_{jk}=-W_{kj}$.
\end{center}

Let $\mathcal{N}=N(N+1)/2$. Any  gravitational configuration 
 $$\mathcal{Q}=(z_1,z_2,\cdots,z_N,w_1,w_2,\cdots,w_N,\delta_{12},\delta_{13},\cdots,\delta_{(N-1)N}), $$
 associates with  two vectors in $\mathbb{C}^{\mathcal{N}}$
 \begin{center}
 	$\mathcal{Z}=(z_1,z_2,\cdots,z_N,Z_{12},Z_{13},\cdots,Z_{(N-1)N})$, $\mathcal{W}=(w_1,w_2,\cdots,w_N,W_{12},W_{13},\cdots,W_{(N-1)N}).$
 \end{center}
  Let $\|\mathcal{Z}\|=\max_{j=1,2,\cdots,\mathcal{N}}|\mathcal{Z}_{j}|$ be the modulus of the maximal component of
the vector $\mathcal{Z}\in \mathbb{C}^\mathcal{N}$. Similarly, set  $\|\mathcal{W}\|=\max_{k=1,2,\cdots,\mathcal{N}}|\mathcal{W}_{k}|$.

Consider
a sequence $\mathcal{Q}^{(n)}$, $n=1,2,\cdots$, of solutions  of (\ref{equ:sys-3}), which are of course normalized central configurations. Take a sub-sequence such that the maximal component of $\mathcal{Z}^{(n)}$ is fixed, i.e., there is a $j\in \{1,2,\cdots,\mathcal{N}\}$ that is  independent 
of $n$ such that  $\|\mathcal{Z}^{(n)}\|=|\mathcal{Z}^{(n)}_{j}|$. 
Extract again in such a way that the sequence $\mathcal{Z}^{(n)}/\|\mathcal{Z}^{(n)}\|$ converges.
Extract again  in such a way that  the maximal component of $\mathcal{W}^{(n)}$ is fixed. Finally,  extract  in
such a way that the sequence $\mathcal{W}^{(n)}/\|\mathcal{W}^{(n)}\|$ converges.

If the initial sequence has the property  that $\mathcal{Z}^{(n)}$ or $\mathcal{W}^{(n)}$ is unbounded, so is the extracted
sequence. Note that $\|\mathcal{Z}^{(n)}\|$ and $\|\mathcal{W}^{(n)}\|$ are bounded away from zero: if the first $N$
components of the vector $\mathcal{Z}^{(n)}$ or $\mathcal{W}^{(n)}$ all approach  zero, then the denominator $z_{12}=w_{12}$ of the
component $Z_{12}=W_{12}$  approaches  zero, hence  $\mathcal{Z}^{(n)}$ and $\mathcal{W}^{(n)}$ are unbounded. There are two possibilities for the
extracted sub-sequences above:\begin{itemize}
	\item $\mathcal{Z}^{(n)}$ and $\mathcal{W}^{(n)}$ are bounded,
	\item at least one of $\mathcal{Z}^{(n)}$ and $\mathcal{W}^{(n)}$ is unbounded.
\end{itemize}
\begin{definition}[Singular sequence]
	Consider a sequence of normalized central configurations. A
	sub-sequence extracted by the above process, in the unbounded case, is called
	a \emph{singular sequence}.
\end{definition}

\section{Rules of colored diagram} \label{sec:rule}

\indent\par

In this section, we state the rules of the colored diagram, which will be the tool for the classification of the singular sequences. Most rules are the same as that in \cite{Albouy2012Finiteness}. We remark that 
 it is presumed that the sum, and partial sum of the masses are nonzero in \cite{Albouy2012Finiteness} and that here  this hypothesis  is discarded. As a result, there is no center of mass. Hence, Rule 1c, 1d, and 1e of \cite{Albouy2012Finiteness} are modified into our Rule III, IV, and V, and Rule 2a  of \cite{Albouy2012Finiteness} has been dropped.

\begin{definition}[Notation of asymptotic estimates]
	\ \ 
\begin{itemize}
\item[$a\sim b$]  means $a/b\rightarrow 1$;
\item[$a\prec b$]   means $a/b\rightarrow 0$;
\item[$a\preceq b$]  means $a/b$ is bounded;
\item[$a\approx b$]  means $a\preceq b$ and $a\succeq b$.
\end{itemize}
\end{definition}

\begin{definition}[Strokes and circles.]
	
Given a singular sequence,  the indices of the bodies  will be written in a figure and we use two colors for edges and
	vertices.
	
	We use the first color, called the $z$-color,  to mark the maximal order components
	of
	\begin{center}
		$\mathcal{Z}=(z_1,z_2,\cdots,z_N,Z_{12},Z_{13},\cdots,Z_{(N-1)N})$.
	\end{center}
	They are the components of
	the converging  vector sequence $\mathcal{Z}^{(n)}/\|\mathcal{Z}^{(n)}\|$ that do not go  to zero. If the term $z^{(n)}_j$ is of maximal order among all the components of $\mathcal{Z}^{(n)}$, we draw a circle around
	the name of vertex $\textbf{j}$. If the term
	$Z^{(n)}_{kl}$ is of maximal order among all the components of $\mathcal{Z}^{(n)}$, we draw a stroke between the names $\textbf{k}$ and $\textbf{l}$.  
\end{definition}

\indent\par

If there
is a maximal order term in an equation, there should be another one.  Hence, we have the first Rule. 
\begin{description}
	\item[{Rule I}]
	At  each end of any $z$-stroke, there is another $z$-stroke
	or/and a $z$-circle drawn around the name of the body. A $z$-circle is not  isolated;  a $z$-stroke must  emanate  from it. There exist  at least one
	$z$-stroke in the $z$-diagram.
\end{description}

\begin{definition}[$z$-close]
	
Given  a singular sequence,  we say that bodies $\textbf{k}$ and $\textbf{l}$
	are close in $z$-coordinate, or $z$-close, or that $z_k$ and $z_l$ are close, if $z^{(n)}_{kl}\prec \|\mathcal{Z}^{(n)}\|$.
\end{definition}

The following statement is obvious.
\begin{description}
	\item[{Rule II}] If bodies $\textbf{k}$ and $\textbf{l}$
	are  close  in $z$-coordinate, they are both $z$-circled or both not
	$z$-circled.
\end{description}

\begin{definition}[Isolated component]\label{isolatedcomponent}
	An isolated component of the $z$-diagram is a subset of vertices
	such that there is no $z$-stroke  between a vertex of this subset and  a vertex of its 
	complement.
\end{definition}

\begin{description}
	\item[{Rule III}]  The moment  of mass  of a set of bodies forming an isolated component of the $z$-diagram is $z$-close to the origin.
\end{description}

\begin{description}
	\item[{Rule IV}]  Consider the $z$-diagram or an isolated component of it. If there
	is a $z$-circle, there is another one. The $z$-circled bodies can not all be
	$z$-close together except that the sum of the masses of these bodies is zero.
\end{description}

\begin{definition}[Maximal $z$-stroke]
	Consider a $z$-stroke between  the names  $\textbf{k}$ and  $\textbf{l}$. It is called 
	a maximal $z$-stroke if $z_k$ and $z_l$ are not close.
\end{definition}

\begin{description}
	\item[{Rule V}]  Consider  a maximal $z$-stroke. At least one of the two ends  is $z$-circled. 
\end{description}

We also draw
$w$-strokes and $w$-circles on the same diagram.  To distinguish from the $z$-strokes and $z$-circles, we use another color. The previous rules
and definitions apply to $w$-strokes and $w$-circles. The superposition of the $z$-diagram and the $w$-diagram will be called  simply  the
diagram.  The definitions about the $z$-diagram will  
be  adapted for the  $w$-digram. For example,  a subset of bodies
is an isolated component of the diagram if and only if it forms an isolated
component of the $z$-diagram and an isolated component of the $w$-diagram.

\begin{definition}[Edges and strokes]
	If there
	is either a $z$-stroke, or a $w$-stroke, or both between vertex $\textbf{k}$ and vertex $\textbf{l}$, we say that there is an edge between them.  There are three types of edges,
	$z$-edges, $w$-edges and $zw$-edges, and  two types of strokes, represented by 
	two different colors.
\end{definition}
\begin{figure}[h!]
	\centering
 \includegraphics[width=12 cm, height= 1.5 cm]{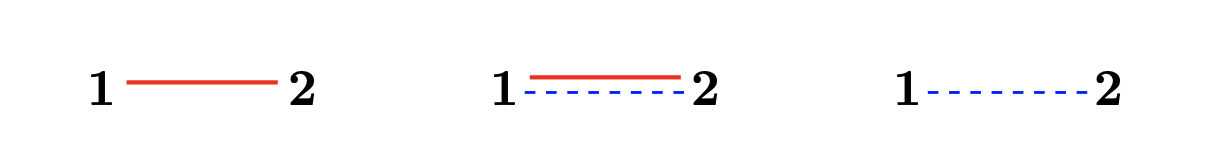}
	\caption{A $z$-stroke, a $z$-stroke plus a $w$-stroke, a $w$-stroke,
		forming respectively a z-edge, a $zw$-edge, a $w$-edge.  }
	\label{fig:edges}
\end{figure}

\subsection{New normalization. Main estimates.}

\indent\par

Note that a central configuration will not change  by multiplying the $z$ coordinates by $a\in\mathbb{C}\backslash\{0\}$ and the $w$ coordinates
by $a^{-1}$.  The  diagram is also invariant under  such an operation, because the $z$-coordinates and the $w$-coordinates are  considered 
separately.

The normalization $z_{12}=w_{12}$ was used  in the previous sections.   In
the following we normalize instead with $\|\mathcal{Z}\|=\|\mathcal{W}\|$. For a central configuration normalized with the condition $z_{12}=w_{12}$, we multiply the $z$-coordinates by a proper $a > 0$, the $w$-coordinates by $a^{-1}$, so that the
maximal component of $\mathcal{Z}$ and the maximal component of $\mathcal{W}$ have the same
modulus, i.e., $\|\mathcal{Z}\|=\|\mathcal{W}\|$.

A singular sequence was defined by the condition either $\|\mathcal{Z}^{(n)}\|$  or $\|\mathcal{W}^{(n)}\|$ tends to $ \infty$.
Recall that both $\|\mathcal{Z}^{(n)}\|$ and $\|\mathcal{W}^{(n)}\|$ were bounded away from zero. With the
new normalization, a singular sequence is simply characterized by $\|\mathcal{Z}^{(n)}\|=\|\mathcal{W}^{(n)}\|\rightarrow \infty$. From now on we
only discuss singular sequences.

Set $\|\mathcal{Z}^{(n)}\|=\|\mathcal{W}^{(n)}\|=1/\epsilon^2$, then $\epsilon\rightarrow 0$.

\begin{proposition}[\cite{Albouy2012Finiteness} Estimate 1]\label{Estimate1}
	For any $(k,l)$, $1\leq k<l\leq N$, we have $\epsilon^2\preceq z_{kl}\preceq \epsilon^{-2}$, $\epsilon^2\preceq w_{kl}\preceq \epsilon^{-2}$ and $\epsilon \preceq r_{kl}\preceq \epsilon^{-2}$.
		There is a $zw$-edge between $\textbf{k}$ and $\textbf{l}$ if and only if $r_{kl}\approx \epsilon$. There is a maximal $z$-edge between $\textbf{k}$ and $\textbf{l}$ if and only if $w_{kl}\approx \epsilon^{2}$.
	\end{proposition}

\begin{proposition}[\cite{Albouy2012Finiteness}  Estimate 2]\label{Estimate2}		
	We assume that there is a $z$-stroke between $\textbf{k}$ and $\textbf{l}$. Then
	\[  \epsilon \preceq r_{kl} \preceq 1, \ \epsilon \preceq z_{kl} \preceq \epsilon^{-2},\    \epsilon \succeq w_{kl} \succeq \epsilon^{2}.\  \]
	Under the same hypothesis the ``equality case'' are characterized as follows:
	\begin{align*}
	{\rm Left:\ } &r_{kl} \approx \epsilon \Leftrightarrow z_{kl} \approx \epsilon \Leftrightarrow w_{kl} \approx \epsilon \Leftrightarrow zw{\rm -edge\ between\ } k \ {\rm and} \ l,\\
		{\rm Right:\ } &r_{kl} \approx 1 \Leftrightarrow z_{kl} \approx \epsilon^{-2} \Leftrightarrow w_{kl} \approx \epsilon^2 \Leftrightarrow {\rm maximal\ } z{\rm -edge\ between\ } k \ {\rm and} \ l.
	\end{align*}
		\end{proposition}

\begin{remark}
	By the estimates above, the strokes in a $zw$-edge are not maximal. A maximal $z$-stroke is exactly a maximal $z$-edge.
\end{remark}

\begin{description}
	\item[{Rule VI}] \emph{Two consecutive $zw$-edges.} If two consecutive $zw$-edges are present, there is a third $zw$-edge closing the triangle.
\end{description}

\emph{Clusters}. At the limit when following a singular sequence, the $z_k$'s form clusters. For example, if bodies 1, 2 and 3 are such that $z_{12}\prec z_{23}$, we say that 1 clusters with 2 in $z$-coordinate, relatively to the bodies 1,2,3. If there is a fourth body such that $z_{24}\prec z_{12}\prec z_{23}$, we say the fourth body  form a sub-cluster, e.g., together with body 2. 

We can  write a \emph{clustering scheme} in each coordinate. For instance, the situation considered above is  simply  $z: 42. 1...3$, where three dots means  the largest separation within the group, one dot  means the intermediate separation, and no dot means  the smallest separation.

In the following rule   clustering relation inside a sub-system of three bodies is considered. Note  that these three bodies may  form, e.g., in $z$-coordinate, a cluster relatively to the whole configuration.

\begin{description}
	\item[{Rule VII}] \emph{Skew clustering.} Suppose that there is an edge from vertex 1 to vertex 2,  an edge from vertex 2 to vertex 3, and that there is no edge from vertex 1 to vertex 3.  Then the clustering schemes are $z: 1.2...3, w: 1...2.3$, or $z: 1...2.3, w: 1.2...3$. We say there is ``skew clustering''.
\end{description}

\begin{corollary} [\cite{Albouy2012Finiteness}]
	Two consecutive $z$-edges cannot be maximal if they are not part of a triangle of edges.
\end{corollary}

\begin{description}
	\item[{Rule VIII}] \emph{Cycles.} Consider  a cycle of edges, the list of $z$-separations corresponding to the edges, and the maximal oder of the $z$-separations within this list. Two or more of the $z$-separations are of this order. The corresponding edges have the same type. If there are only two, the corresponding separations are not only of the same order, but equivalent.
\end{description}

\begin{description}
	\item[{Rule IX}] \emph{Triangles.}  Consider a triangle of edges in the diagram. Then the edges are of  the same type (all $z$-edges or all $w$-edges or all $zw$-edges), all the $z$-separations are of the same order, all the $w$-separations are of the same order.
\end{description}

\begin{corollary} [\cite{Albouy2012Finiteness}]
	Consider three vertices. There are 6, 3, 2, 1, or 0 strokes joining them. If there are three forming a triangle, they are of the same color.
\end{corollary}

\begin{description}
	\item[{Rule X}] \emph{Fully edged sub-diagram.}  Suppose that in the diagram there is a triangle of edges, plus a fourth vertex attached to the triangle by at least two edges, plus a fifth vertex attached to the four previous vertices by at least two edges, and so on up to a $p$-th vertex, $p\ge 3$. Then there is indeed an edge, of the same type,  between any pair of the $p$ vertices, all the $z$-separations are of the same order, all the $w$-separations are of the same order.
\end{description}

\begin{description}
	\item[{Rule XI}] If four edges form a quadrilateral, then the opposite edges are of the same type.
\end{description}

\begin{description}
	\item[{Rule XII}] \emph{Bounded potential.}  Consider a singular sequence. Let  $k_0$ and $l_0$ be such a pair of  bodies  that $r_{k_0l_0} \preceq r_{kl}$ for any $k, l, 1\le k<l\le n$. If $r_{k_0l_0}\longrightarrow0$, then there is another pair of bodies $(k_1, l_1)$ such that  $r_{k_0l_0}\approx r_{k_1l_1}$.
\end{description}

\begin{corollary}\label{cor:zw-edge-not-alone}
	If a  $zw$-edge is present  in the diagram, there is another one.
\end{corollary}

\section{ Exclusion of 4-body diagrams}\label{sec:list}
\indent\par
A \emph{bicolored vertex} of the diagram  is a vertex which connects at least
a $z$-stroke  with at least a $w$-stroke. The number of edges from
a bicolored vertex is at least 1 and at most $3$. The number of strokes from
a bicolored vertex is at least 2 and at most $6$. Given a diagram, we
define $C$ to be  the maximal number of strokes from a bicolored vertex. 
This number will be used  to classify all possible diagrams.

Recall that the $z$-diagram indicates the maximal order terms and it is nonempty. If there is a circle, there is an edge of the same
color adjacent to the circle. So there is at least a $z$-stroke, and at least a
$w$-stroke.

\subsection{No bicolored vertex}

\indent\par

\begin{figure}[h!]
	\centering
	 \includegraphics[width=5 cm, height= 3.3 cm]{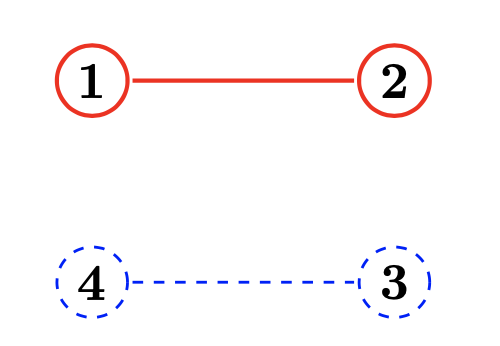}
	\caption{No bicolored vertex}
		\label{fig:C=0}
\end{figure}

If there does not exist a  bicolored vertex, then
there are at most two strokes and they are ``parallel''.  Thus
the only possible diagram is the one in  Figure \ref{fig:C=0}.

\subsection{$C=2$}

\indent\par
There are two cases: a $zw$-edge exists or not.

\begin{figure}[h!]
	\centering
	 \includegraphics[width=4.6 cm, height= 3.3 cm]{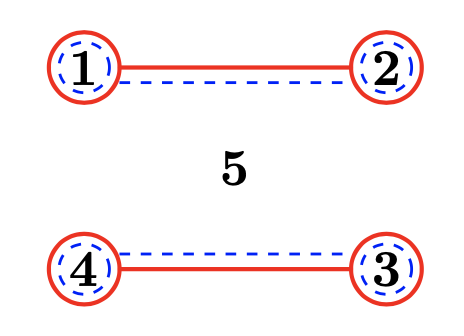}

	\caption{$C=2$, $zw$-edge appears }
	\label{fig:C=21}
\end{figure}

If it is present,  it should be isolated.  On the other hand, there should be another $zw$-edge by Corollary \ref{cor:zw-edge-not-alone}. Then
the only possible diagram  is the one  in  Figure \ref{fig:C=21}.  By Rule IV, we have $m_1+m_2=0$ since body 1 and 2 are both $z$-circled and $z$-close. Similarly, we get  $m_3+m_4=0$, then $\sum_{i=1}^4m_i=0$.

\begin{figure}[h!]
	\centering
	 \includegraphics[width=4.5 cm, height= 3 cm]{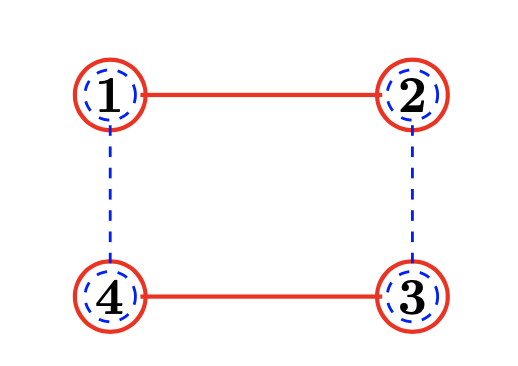}
	
	\caption{$C=2$, no $zw$-edge  }
	\label{fig:C=22}
\end{figure}

If it is not present, there are adjacent $z$-edges and $w$-edges. From any
such adjacency there is no other edge.  By trying to continue it, we see
that the only diagram is the one in  Figure \ref{fig:C=22}.

\subsection{$C=3$}

\indent\par
Consider a bicolored vertex  with three strokes.
There are two cases: a $zw$-edge exists or not.

If it is not present, it is Y-shaped. Suppose that vertex $\textbf{1}$ connects with vertex $\textbf{2}$ and vertex $\textbf{3}$ by   $z$-edges, and  connects with  vertex $\textbf{4}$ by  a $w$-edge. By Rule I, vertex $\textbf{1}$ is $w$-circled, then vertex $\textbf{2}$ and vertex $\textbf{3}$ are also $w$-circled by Estimate 2. By Rule I again, there is a   $w$-stroke emanating from vertex $\textbf{2}$ and vertex $\textbf{3}$, which leads to a triangle with edges of different types. We exclude this case by Rule IX.

\begin{figure}[h!]
		 \includegraphics[width=15 cm, height= 3 cm]{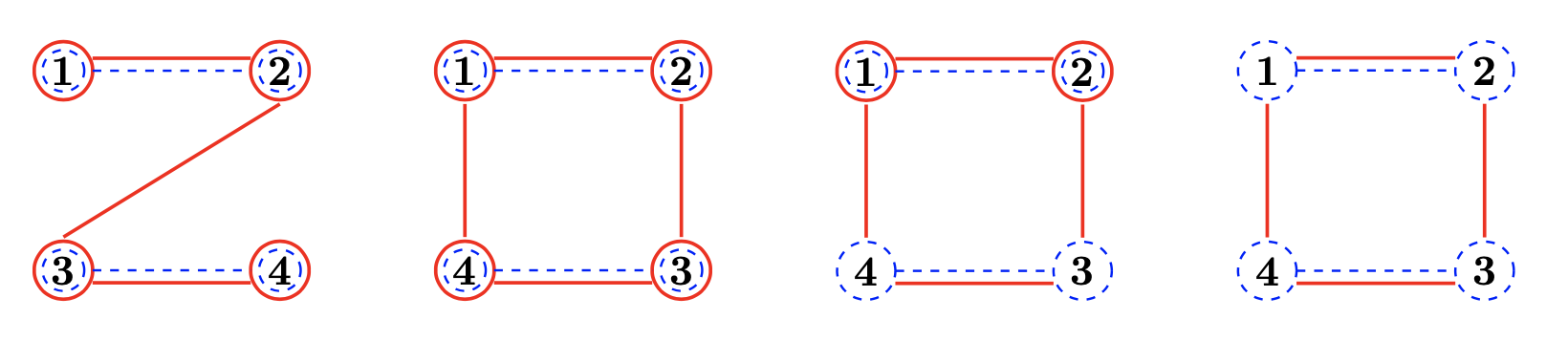}
	\caption{$C=3$}
	\label{fig:C=3}
\end{figure}

If it is present, let vertex $\textbf{2}$ be the bicolored vertex with three strokes. Suppose it connects with vertex $\textbf{1}$ by  a $zw$-edge,  with vertex $\textbf{3}$ by a  $z$-edge. By circling method, we circle the three vertices by $w$-color. Then there is $w$-stroke from vertex $\textbf{3}$, which can not connect to $\textbf{1}$ by Rule IX. Then there are two cases: an edge between vertex $\textbf{1}$ and vertex $\textbf{4}$ or not. If it is not present, then we continue it to the first one of  Figure \ref{fig:C=3}. If it is present, then they form a quadrilateral. By Rule XI, we see the types of edges.  There is no diagonal  by
Rule IX. Then circling method gives the next three diagrams  in  Figure \ref{fig:C=3}.  Consider the isolated components of the $w$-diagram. We see  $m_1+m_2=0$ and $m_3+m_4=0$ by Rule IV, for all the four diagrams in  Figure \ref{fig:C=3}.

\subsection{$C=4$}

\indent\par
Consider a bicolored vertex (let us say, $\textbf{1}$ ) with four strokes.

\begin{figure}[h!]
	\centering
	 \includegraphics[width=12.5 cm, height= 3.2 cm]{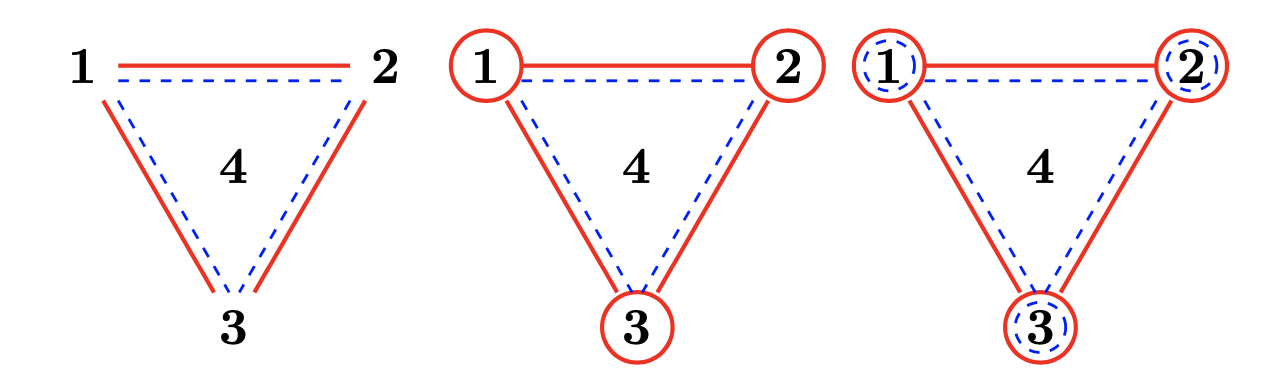}
	
	\caption{$C=4$, three $zw$-edges}
	\label{fig:C=41}
\end{figure}

In the first case, vertex $\textbf{1}$  has  $zw$-edges connected with  vertex $\textbf{2}$ and  vertex $\textbf{3}$ separately. A third
$zw$-edge closes the triangles by Rule VI.  As $C=4$,  there is no other  stroke, thus vertex  $\textbf{4}$ is neither $z$-circled nor $w$-circled.  Hence, 
the possible diagrams are those  in Figure \ref{fig:C=41}.

In the second case, vertex $\textbf{1}$  has one adjacent $zw$-edge connected with  vertex $\textbf{2}$, a $z$-edges connected with vertex $\textbf{3}$ and a $w$-edges connected  with vertex $\textbf{4}$. Any other edge in this diagram would close a triangle, which would contradicts Rule IX. Hence, $\textbf{1}$ is $z$ and $w$-circled, and so is vertex $\textbf{3}$ by Rule I and Rule II. Then there is a  $w$-edge emanating from $\textbf{3}$ by Rule I,  contradiction.

In the third case, vertex $\textbf{1}$  has one adjacent $zw$-edge connected with  vertex $\textbf{2}$, two $z$-edges connected with vertex $\textbf{3}$ and $\textbf{4}$. If there are more strokes, it should be a $z$-stroke between
vertex $\textbf{3}$ and $\textbf{4}$ by Rule IX. However, $\textbf{2}$ is $z$ and $w$-circled, and so is vertex $\textbf{1}$ by Rule II. Then both of vertex $\textbf{3}$ and $\textbf{4}$ are $w$-circled by Rule II. Then there should be  a  $w$-stroke emanating from $\textbf{3}$, contradiction.

\subsection{$C=5$}

\indent\par
Consider a bicolored vertex (let us say, $\textbf{1}$ ) with five strokes.

Suppose  vertex $\textbf{1}$  has one adjacent $w$-edge connected with   vertex $\textbf{4}$, and two adjacent $zw$-edges connected with  vertex $\textbf{2}$ and  vertex $\textbf{3}$. Then there is  a fully $zw$-edged triangle between vertexes  $\textbf{1,2,3}$ by Rule VI.  Rule VI  also implies that there are no more edges. By Rule I and IV, all vertices are $w$-circled. There is no $z$-circle, otherwise all
vertices are $z$-circled by Rule IV and II. Then there would be a  $z$-edge emanating from vertex $\textbf{4}$, a contradiction.

\begin{figure}[h!]
	\centering
 \includegraphics[width=5.5 cm, height= 3.5 cm]{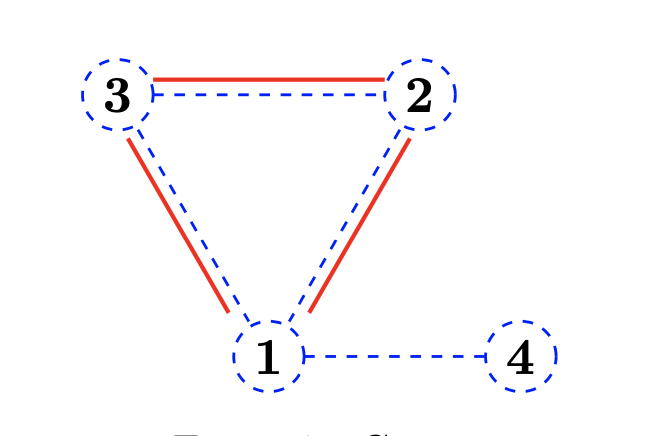}
	\caption{$C=5$}
		\label{fig:C=5}
\end{figure}

\subsection{$C=6$}

\indent\par
\begin{figure}[h!]
	\centering
 \includegraphics[width=13 cm, height= 3.6 cm]{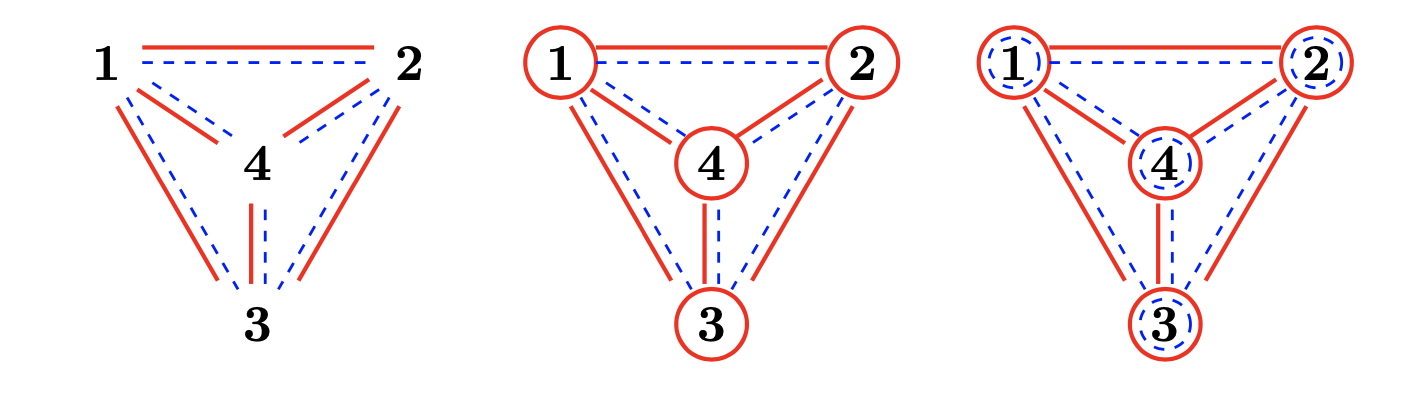}
	\caption{$C=6$}
	\label{fig:C=61}
\end{figure}

Consider a bicolored vertex  with six strokes. Then   this is  a fully $zw$-edged diagram by Rule VI. According to if there are $z$-circle or $w$-circle at vertices, the  possible diagrams are those   in  Figure \ref{fig:C=61}.

\section{Problematic diagrams with non-vanishing total mass} \label{sec:dia}

\emph{Notations.} From now on, $x^{\frac{1}{n}}$ (or $\sqrt[n]{x}$) is understood as  one appropriate value  of the $n$-th root of $x$. 
We denote the set $f=0$ by $\mathcal{V}_f$, and $f=0,g=0$ by $\mathcal{V}_f\bigcap \mathcal{V}_g$ or simply $\mathcal{V}_f \mathcal{V}_g$; the set   $f=0$ or $g=0$ (i.e., $fg=0$) by $\mathcal{V}_f\bigcup \mathcal{V}_g$. Recall that $m$ is the total mass and it is nonzero. 

It is easy to see that all   diagrams  in  Figure  \ref{fig:C=21} and  Figure  \ref{fig:C=3} and the last two diagrams  in  Figure \ref{fig:C=61} must have  $m=0$.  Now, We show that the second diagram in Figure \ref{fig:C=41} is also impossible.

First we have
\begin{equation}\label{d0}
m_1+m_2+m_3=0. 
\end{equation}

The equations $w_k=\sum_{j\neq k} m_j W_{jk}, k=1,2,3,$ imply 
\begin{equation}\label{d42}
\   \frac{W_{12}}{m_3}    \sim \frac{W_{23}}{m_1}  \sim \frac{W_{31}}{m_2}\sim a \epsilon^{-2},  \ \  a\ne 0. 
\end{equation}

Note that 
\[  -m_4 w_4 = \sum_{j=1}^{3} m_j w_j = m_2 w_{12} + m_3 w_{13}=m_1w_{21}+m_3 w_{23} = m_1w_{31} +m_2 w_{32}\preceq \epsilon.  \]
We claim that $w_4 \approx \epsilon$. Otherwise, the above equation leads to 
\[  \frac{w_{12}}{m_3} \sim  \frac{w_{23}}{m_1} \sim   \frac{w_{31}}{m_2} \sim   b \epsilon,  \ \  b\ne 0.   \]
Then equation \eqref{d42} and the identity $W_{kl}= w_{kl}^{-1/2} z_{kl}^{-3/2}$ 
imply that  
\[ z_{12}\sim (\frac{1}{m_3^3 a^2 b})^{1/3} \epsilon, \ z_{23}\sim (\frac{1}{m_1^3 a^2 b})^{1/3} \epsilon, z_{31}\sim (\frac{1}{m_2^3 a^2 b})^{1/3} \epsilon.   \]
Then the identity $z_{12}+z_{23}+z_{31}=0$ implies that 
\[  \frac{1}{m_1} 1^{\frac{1}{3}} + \frac{1}{m_2} 1^{\frac{1}{3}} +\frac{1}{m_3} 1^{\frac{1}{3}} =0,   \]
which contradicts with equation \eqref{d0}. 

Note that 
\[ \sum_{j=1}^3 m_j w_{j4}  =m w_4  \approx \epsilon, \  w_{14}, w_{24}, w_{34} \prec \epsilon^{-2}, \ z_{14} \approx \epsilon^{-2}, \ z_{24}- z_{14}\approx z_{34}-z_{14} \approx \epsilon.    \]
Hence, 
\begin{align*}
m I/m_4 &\sim    \sum_{j=1}^3 m_j z_{j4} w_{j4}\\
& =   z_{14}\sum_{j=1}^3 m_j w_{j4} + (z_{24}- z_{14}) m_2 w_{24}+ (z_{34}-z_{14}) m_3 w_{34}\\
&\approx \epsilon^{-1}. 
\end{align*}
This is a contradiction since the momentum of inertia $I$ should be bounded.

\begin{figure}[t!]
	\centering
	 \includegraphics[width=8.5 cm, height= 9.5 cm]{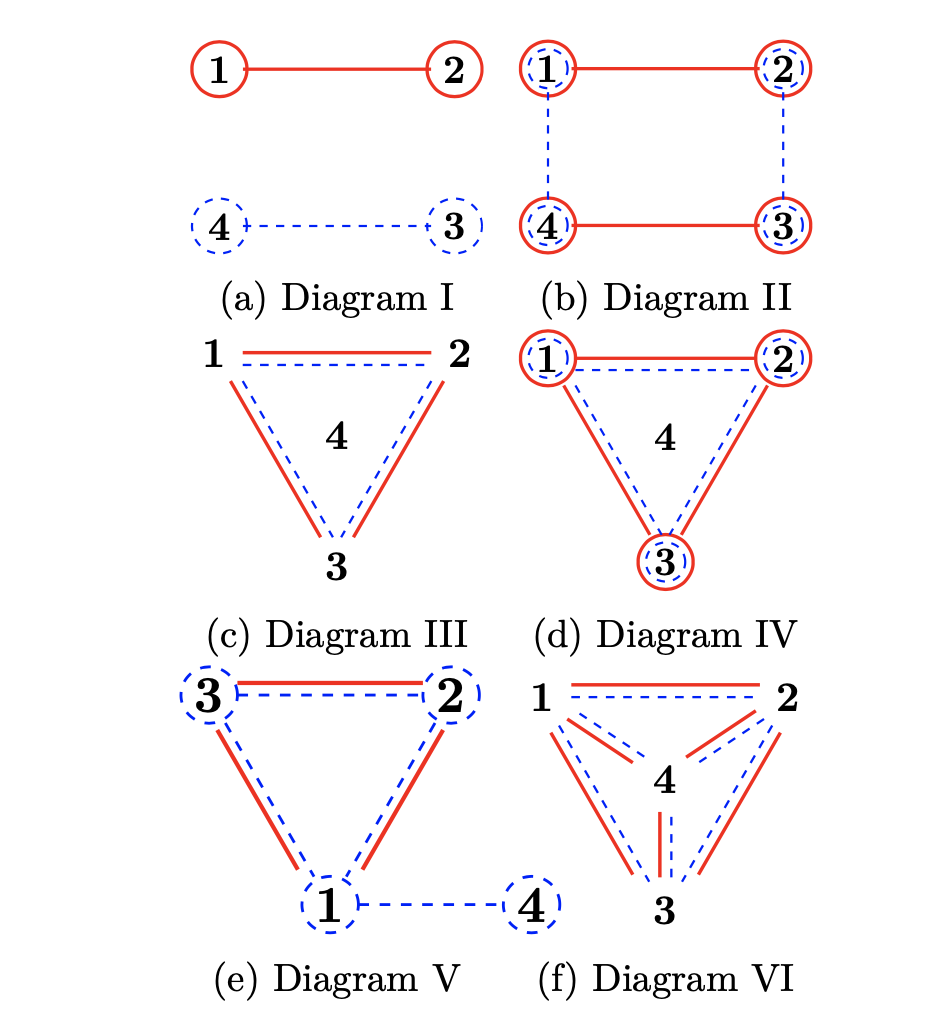}
\caption{Problematic Diagrams}

	\label{fig:Problematicdiagrams}
\end{figure}

We could not exclude  the
diagrams in Figure \ref{fig:Problematicdiagrams}. Some singular sequence could still exist and approach
any of these diagrams.   

\subsection{Diagram I}\label{subsec:DiagramI}

Following the argument from \cite{Albouy2012Finiteness}, we obtain the equation
 \begin{equation} \notag
 \frac{m_1m_3}{(m_2m_4)^{\frac{1}{2}}} +  \frac{m_2m_3}{(-m_1m_4)^{\frac{1}{2}}} +  \frac{m_1m_4}{(-m_2m_3)^{\frac{1}{2}}} +  \frac{m_2m_4}{(m_1m_3)^{\frac{1}{2}} } =0,
\end{equation}
or equivalently,
 \begin{equation} \label{equ:e-D1}
m_1m_3\sqrt{m_1m_3} + m_2m_3\sqrt{-m_2m_3} + m_1m_4\sqrt{-m_1m_4}   + m_2m_4\sqrt{m_2m_4}      =0.
\end{equation}

Note that $z_{jk}, w_{jk}\approx \epsilon^{-2}$, $j=1,2,k=3,4$.  It follows 
\begin{equation}\label{DiagramI4}
r_{jk}\approx \epsilon^{-2},
j=1,2,k=3,4.
\end{equation}
and $Z_{jk}, W_{jk}\approx \epsilon^{4}$. Combined with the following two equations 
\begin{equation}\notag \label{DiagramI2}
    \begin{array}{c}
      w_{12}= (m_1+m_2) W_{12}+m_3(W_{32}-W_{31})+m_4(W_{42}-W_{41}), \\
      z_{34}= (m_3+m_4) Z_{34}+m_1(Z_{14}-Z_{13})+m_4(Z_{24}-Z_{23}), 
    \end{array}
\end{equation}
we obtain 
\begin{equation}\label{DiagramI1}
(m_1+m_2) (m_3+m_4)\neq 0,    r_{12}^3\sim m_1+m_2,
r_{34}^3\sim m_3+m_4.
\end{equation}

Let us  simplify  equation \eqref{equ:e-D1}.  Note that equation \eqref{equ:e-D1} is homogeneous, and is invariant  
 under the transformation
 \[ (m_1,m_2, m_3, m_4)\mapsto (m_2, m_1, m_4, m_3). \]
 Assume that the signs of the masses are the same. It suffices to consider only the case  $(+,+,+,+)$. Then equation  \eqref{equ:e-D1} becomes
$$m_1m_3\sqrt{m_1m_3}  + m_2m_4\sqrt{m_2m_4}  + i (m_2m_3\sqrt{m_2m_3} + m_1m_4\sqrt{m_1m_4}  )=0,$$
which is $m_1 ^3m_3^3 =m_2^3m_4^3$ and $m_2 ^3m_3^3 =m_1^3m_4^3$. Hence, it is  necessary that 
 \begin{equation} \notag
m_1=m_2, \ \ m_3=m_4, \    m_1 m_3 >0.
\end{equation}

Assume that only one of the signs of the masses is different from the others. It suffices to consider only   $(+,+,+,-)$ and $(+,-,+,+)$. In the first subcase,
equation \eqref{equ:e-D1}  is equivalent to
\begin{equation*}\begin{cases}
m_1m_3\sqrt{m_1m_3}  + m_1m_4\sqrt{-m_1m_4} =0, \cr
 m_2m_3\sqrt{-m_2m_3}  + m_2m_4\sqrt{m_2m_4} =0. 
\end{cases}
\end{equation*}
Then we have $m_3=-m_4$, which contradicts with equation \eqref{DiagramI1}. Similarly,  the second subcase is impossible neither. 

 Assume that only two  of the signs of the masses is different from the others. It suffices to consider only   $(+,+,-,-)$, $(+,-,+,-)$ and $(+,-,-,+)$. In the first subcase,
 equation \eqref{equ:e-D1}  is equivalent to
 \begin{equation*}\begin{cases}
m_1m_3\sqrt{m_1m_3}  + m_2m_4\sqrt{m_2m_4} =0, \cr
m_2m_3\sqrt{-m_2m_3}  + m_1m_4\sqrt{-m_1m_4} =0,
 \end{cases}
 \end{equation*}
  which is  equivalent to
 \begin{equation} \notag
m_1=m_2, \ m_3=m_4, \ m_1 m_3 <0.
\end{equation}
In the rest subcases, equation \eqref{equ:e-D1}  can not be simplified.  

 To summarize, the masses of central configurations corresponding to the first diagram belong to one of the following two sets of constraints
\begin{center}
$\begin{array}{ll}
    \mathcal{V}_{IA}[12,34]: & \hbox{$m_1=m_2, \ m_3=m_4;$} \\
  \mathcal{V}_{IB}[12,34]: &\hbox{$ m_1m_3\sqrt{m_1m_3} + m_2m_3\sqrt{-m_2m_3} + m_1m_4\sqrt{-m_1m_4}   + m_2m_4\sqrt{m_2m_4}      =0,$}\\
 &\hbox{$ m_1m_2<0, \ m_3m_4<0.$}
  \end{array}
$
\end{center}

  Furthermore, we claim that 
  if $I=0$, then the masses belong to the set
  \begin{equation}\label{DiagramIIis0}
    \mathcal{I}_{I}[12,34]:~~~~~~~~~ \frac{m_1m_2}{\sqrt[3]{m_1+m_2}}+\frac{m_3m_4}{\sqrt[3]{m_3+m_4}}=0. 
  \end{equation}
   This is a consequence of the fact that $U=I=0$.
   
In fact, it is easy to see that if $I=0$, then the masses belong to the set \begin{equation}\label{DiagramIIis0new} \notag
   \mathcal{V}_{I0}[12,34]~~~~\triangleq~~~~~ \mathcal{V}_{IB}[12,34] \mathcal{I}_{I}[12,34].
  \end{equation}

Note that Diagram I has five other variants where the vertices are renumbered. Two of them are presented in Figure \ref{fig:sec6.1}. We still call them Diagram I if there is no confusion arise.  The corresponding mass sets are denoted by 
\[ \mathcal{V}_{IA}[ij,pk], \   \mathcal{V}_{IB}[ij,pk], \  \mathcal{I}_{I}[ij,pk],   \mbox{ and }   \mathcal{V}_{I0}[ij,pk],  \]
  if there is one stroke, either $z$ or $w$-stroke,  between vertices \textbf{i} and \textbf{j}, and the other stroke is between  
  vertices \textbf{p} and \textbf{k}. For instance,  for both of the two  diagrams in Figure \ref{fig:sec6.1}, the corresponding mass sets are denoted by $\mathcal{V}_{IA}[13,24],    \mathcal{V}_{IB}[13,24],  \mathcal{I}_{I}[13,24],$ or $ \mathcal{V}_{I0}[13,24]$.

\subsection{Diagram II}
By the argument from \cite{Albouy2012Finiteness} and the   hypothesis $m\ne 0$ , the masses of central configurations corresponding to the second diagram belong to the following set
\begin{align*}
 \mathcal{V}_{II}[13,24]:~~~~~~~~~&m_1m_3=m_2m_4, \\ &(m_1+m_2) (m_2+m_3)  (m_3+m_4)(m_1+m_4)\neq 0.
\end{align*}

By the facts $ z_{12}= (m_1+m_2) Z_{12}+m_3(Z_{32}-Z_{31})+m_4(Z_{42}-Z_{41})$ and $Z_{12}=z_{12}r_{12}^{-3}$,  it follows  $r_{12}^3\sim (m_1+m_2)$. Similarly, we have 
\begin{equation}\label{DiagramII2}
 r_{12}^3\sim (m_1+m_2),  r_{23}^3\sim (m_2+m_3), r_{34}^3\sim (m_3+m_4),  r_{14}^3\sim (m_1+m_4).
\end{equation}
Note that   \begin{equation}\label{DiagramII3}
   r_{13}\approx \epsilon^{-2},
    r_{24}\approx \epsilon^{-2}.
\end{equation}

If $I=0$, then the masses belong further to the set
  \begin{equation}\label{DiagramIIIis0}
     \mathcal{I}_{II}[13,24]:~~~~~~~~~  \frac{m_1m_2}{\sqrt[3]{m_1+m_2}}+\frac{m_2m_3}{\sqrt[3]{m_2+m_3}}+\frac{m_3m_4}{\sqrt[3]{m_3+m_4}}+\frac{m_1m_4}{\sqrt[3]{m_1+m_4}}=0. 
  \end{equation}
  That is, if $I=0$, then the masses belong to the set \begin{equation}\label{DiagramIIIis0new}\notag
   \mathcal{V}_{II0}[13,24]~~~~\triangleq~~~~~  \mathcal{V}_{II}[13,24] \mathcal{I}_{II}[13,24].
  \end{equation}
This is a consequence of the fact that $U=I=0$.

Similarly,  Diagram II  has five other variants. Some of them are presented in the second row of Figure \ref{fig:Problematicdiagramsr21cj}. We still call them Diagram II if there is no confusion arise.  The corresponding mass sets are denoted by 
\[  \mathcal{V}_{II}[ij,pk], \  \mathcal{I}_{II}[ij,pk],     \mbox{ and }    \mathcal{V}_{II0}[ij,pk],  \]
if one  diagonal is  between vertices \textbf{i} and \textbf{j}, and the other  is between  
vertices \textbf{p} and \textbf{k}. For instance,  for both of the last  two  diagrams in  the second row of Figure \ref{fig:Problematicdiagramsr21cj},  the corresponding mass sets are denoted by $ \mathcal{V}_{II}[14,23], \mathcal{I}_{II}[14,23]$ or $ \mathcal{V}_{II0}[14,23]$.

\subsection{Diagram III}

Following the argument from \cite{Albouy2012Finiteness}, we arrive at the estimation
\[     \frac{Z_{12}}{m_3}    \sim \frac{Z_{23}}{m_1}  \sim \frac{Z_{31}}{m_1}, \   \frac{W_{12}}{m_3}    \sim \frac{W_{23}}{m_1}  \sim \frac{W_{31}}{m_1}.  \]
Since $Z_{12}, W_{12}$ are both of the maximal order, we assume the first set is asymptotic to $a \epsilon^{-2}$ and the second asymptotic to $b \epsilon^{-2}$, $ab\ne 0$. Multiplying the two, by $Z_{kl}W_{kl}= r_{kl}^{-4}$, we obtain
\[ \frac{1}{m_3^2 r_{12}^4} \sim   \frac{1}{m_1^2 r_{23}^4} \sim \frac{1}{m_2^2 r_{31}^4} \sim ab \epsilon^{-4}.    \]
Then we have
\[r_{12}\sim 1^{\frac{1}{4}}\frac{\epsilon}{ \sqrt{|m_3|}}c, \ r_{23}\sim 1^{\frac{1}{4}} \frac{\epsilon}{ \sqrt{|m_1|}}c, \ r_{31}\sim 1^{\frac{1}{4}} \frac{\epsilon}{ \sqrt{|m_2|}} c, \ c\ne 0. \]
Then by $z_{kl}=r_{kl}^3 Z_{kl}$, we have
\begin{equation}\label{III-1}
z_{12}\sim 1^{\frac{1}{4}}\frac{\epsilon}{ \sqrt{|m_3|}}d, \  z_{23}\sim 1^{\frac{1}{4}} \frac{\epsilon}{ \sqrt{|m_1|}}d, \ z_{31}\sim 1^{\frac{1}{4}} \frac{\epsilon}{ \sqrt{|m_2|}} d, \ d\ne 0. 
\end{equation}
 Then the identity
$z_{12}+z_{23}+z_{31}=0$ implies that  the masses  belong to the set
\begin{equation} \label{equ:m-D3}  \notag \mathcal{V}_{III}[123]:~~~~~~~~~ \frac{1}{ {|m_1|^{\frac{1}{2}}}} + \frac{1}{ {|m_2|^{\frac{1}{2}}}} +\frac{1}{ {|m_3|^{\frac{1}{2}}}}  =0.   \end{equation}

It is easy to see that 
 \begin{equation}\label{DiagramIII3}
   r_{12},r_{23},r_{13}\approx \epsilon,
    r_{14},r_{24},r_{34}\succ \epsilon.
\end{equation}
We discuss the distances $\{r_{j4}\}_{j=1}^3$ further in the following.  
Without loss of generality, assume that 
\[ z_{14} \preceq z_{24}, z_{34}, \ {\it and}\ z_{14} \preceq w_{14}.   \]
Then $\epsilon \prec w_{14}\sim w_{24}\sim w_{34}$.  For $z_{14}$, there  are three possibilities, namely $\prec, \approx$, or $ \succ \epsilon$. Note that 
\begin{align*}
&z_{14} \prec \epsilon, \ \Rightarrow  z_{24}\approx z_{34}\approx\epsilon, \ r_{14}\prec r_{24}\approx r_{34}, \ W_{14}\succ W_{24}\approx W_{34};\\
&z_{14} \approx\epsilon, \ \Rightarrow  z_{14} \approx z_{24}\approx z_{34}\approx\epsilon, \ r_{14}^2\approx r_{24}^2\approx r_{34}^2, \ W_{14}\approx W_{24}\approx W_{34};\\
&z_{14} \succ \epsilon, \ \Rightarrow  z_{14} \sim z_{24}\sim z_{34}\succ\epsilon, \ r_{14}^2\sim r_{24}^2\sim r_{34}^2, \ W_{14}^2\sim  W_{24}^2\sim W_{34}^2. 
\end{align*}

\subsubsection{The masses for  $z_{14}\prec \epsilon$}

By $w_4=\sum_{j=1}^3 m_jW_{j4}$ and $z_4=\sum_{j=1}^3 m_jZ_{j4}$ it follows that
$
  w_4\sim  m_1W_{14},z_4\sim  m_1 Z_{14}.
$
On the other hand, note that  $\sum_{j=1}^4 m_j w_4=\sum_{j=1}^3 m_jw_{j4}\preceq w_{14}$, then $w_4 \sim m_1 W_{14} \preceq w_{14}$, which implies that $r_{14}\succeq 1$. Therefore $z_4 \sim m_1 Z_{14}\preceq z_{14}$, and 
\[  z_4, z_1\prec \epsilon \approx z_2\approx z_3.  \] 
Then it is easy to see that $ z_{12}\sim z_2\sim -m_3 f \epsilon, z_{13}\sim z_3\sim m_2 f \epsilon$  and $ z_{23}=z_3- z_2\sim (m_2+m_3) f\epsilon$, where $f$ is some nonzero constant. Compared with equation \eqref{III-1},  it is necessary that 
\begin{equation} 
\label{equ:m-D3a}  
4|m_1| =|m_2|=|m_3|. 
\end{equation}

\subsubsection{ When $\sum_{j=1}^3 m_j\neq 0$}
The identity  $\sum_{j=1}^4 m_j w_4=\sum_{j=1}^3 m_jw_{j4}$ implies  that 
\begin{equation}\notag
   w_{4}\sim \frac{\sum_{j=1}^3 m_j}{\sum_{j=1}^4 m_j}w_{14}\succ \epsilon.
\end{equation}
By the equation $w_4=\sum_{j=1}^3 m_jW_{j4}$,  
 it follows that
 
\textbf{Case 1).}  If  $ z_{14}\prec \epsilon$,  then $w_4 \sim m_1 W_{14} \sim \sum_{j=1}^3m_j w_{14}$. Then $r_{14}\approx 1$. Thus 
\[   r_{14}^2 \approx 1 \prec r^2_{24} \approx r^2_{34} \prec \epsilon^{-1},  \  4|m_1| =|m_2|=|m_3|; \]

\textbf{Case 2).}  If  $ z_{14}\approx\epsilon$,  then $w_{14}\preceq W_{14}$. Thus 
\[  \epsilon^2 \prec r^2_{14}\approx r^2_{24}\approx r^2_{34}\preceq 1;  \]

\textbf{Case 3).} If $z_{14} \succ \epsilon$,  then $w_{14}\preceq W_{14}$. Thus 
\[  \epsilon^2 \prec r^2_{14}\approx r^2_{24}\approx r^2_{34}\preceq 1.  \]

\subsubsection{When $\sum_{j=1}^3 m_j= 0$}\label{subsec:III-3}
Then it is necessary that 
\begin{equation}
\label{equ:m-D3b}  
\frac{1}{ {|m_1|^{\frac{1}{2}}}} + \frac{1}{ {|m_2|^{\frac{1}{2}}}} +\frac{1}{ {|m_3|^{\frac{1}{2}}}}  =m_1+m_2+m_3=0. 
\end{equation}
The above system has six solutions  of $(m_1, m_2, m_3)$. It is not necessary to find them  for our purpose. 
By $\sum_{j=1}^4 m_j w_4=\sum_{j=1}^3 m_jw_{j4}$ it follows that
$w_{4}\prec w_{14}. $

\textbf{Case 1).}  If  $ z_{14}\prec \epsilon$,  this is impossible since the system by  equations \eqref{equ:m-D3a}  and  \eqref{equ:m-D3b} has no real solution.  

\textbf{Case 2).}  If  $ z_{14}\approx\epsilon$,  then 
\[  \epsilon^2 \prec r^2_{14}\approx r^2_{24}\approx r^2_{34}\prec \epsilon^{-1};  \]

\textbf{Case 3).} If $z_{14} \succ \epsilon$,  then we have 
\[  \epsilon^2 \prec r^2_{14}\approx r^2_{24}\approx r^2_{34}\prec \epsilon^{-4}.  \]

Similarly,  Diagram III  has three other variants. Two of them are presented in the third row of Figure \ref{fig:Problematicdiagramsr21cj}. We still call them Diagram III if there is no confusion arise.  The corresponding mass sets are denoted by 
\[  \mathcal{V}_{III}[ijk], \]
if  the fully edged triangle has 
 vertices \textbf{i},  \textbf{j}, and  \textbf{k}. For instance,  for the   two  diagrams in  the third row of Figure \ref{fig:Problematicdiagramsr21cj},  the corresponding mass sets are denoted by $ \mathcal{V}_{III}[234]$ and $\mathcal{V}_{III}[134]$ respectively. To distinguish, we may refer to the first one as  Diagram III with $\triangle_{234}$,  the second one as  Diagram III with $\triangle_{134}$ and so on.

\subsection{Diagram IV}
The masses  belong to the set
\[  \mathcal{V}_{IV}[123]:~~~~~~~~~  \sum_{j=1}^3 m_j=0,\]
and we have 
\begin{equation}\label{DiagramIV0}
r_{12},r_{23},r_{13}\approx \epsilon, r_{14}\sim \pm r_{24}\sim \pm r_{34}\approx \epsilon^{-2}. 
\end{equation}

Similarly,  Diagram IV  has three other variants. One of them is  presented in the fourth row of Figure \ref{fig:Problematicdiagramsr21cj}. We still call them Diagram IV if there is no confusion arise.  The corresponding mass sets are denoted by 
\[  \mathcal{V}_{IV}[ijk], \]
if  the fully edged triangle has 
vertices \textbf{i},  \textbf{j}, and  \textbf{k}. For instance,  for the   second  diagram in  the fourth row of Figure \ref{fig:Problematicdiagramsr21cj},  the corresponding mass sets is denoted by $ \mathcal{V}_{IV}[124]$.  To distinguish, we may refer to the first one as  Diagram IV with $\triangle_{123}$ and so on.

\subsection{Diagram V}
Following the argument from \cite{Albouy2012Finiteness}, we have 
\[     \frac{Z_{12}}{m_3}    \sim \frac{Z_{23}}{m_1}  \sim \frac{Z_{31}}{m_1}\sim a \epsilon^{-2}, \  {\rm and \ set } \  z_{12}=m_3, \ z_{31} =m_2, \]
with $a$ being some nonzero constant. By $Z_{kl}=r_{kl}^{-3} z_{kl}$, we see
\[   \frac{m_1}{r_{12}^3}    \sim \frac{-(m_2+m_3)}{r_{23}^3}  \sim \frac{m_1}{r_{31}^3}\sim am_1 \epsilon^{-2}.  \]
Then we have
$$r_{12}\sim \epsilon^{\frac{2}{3}} c \sqrt[3]{m_1}, \ r_{23}\sim \epsilon^{\frac{2}{3}} c \sqrt[3]{-(m_2+m_3)}, \ r_{31}\sim \epsilon^{\frac{2}{3}} c\sqrt[3]{m_1} , \ c\ne 0. $$
In this subsection and in this subsection only,  $\sqrt[3]{x} \ (x\in \mathbb{R})$ is understood as the real cubic root of $x$. 
Then the asymptotic relation $\frac{1}{m_1r_{23}} +\frac{1}{m_2r_{31}}+\frac{1}{m_3r_{12}}\sim 0$ implies
\[   \frac{1}{m_1 \sqrt[3]{-(m_2+m_3)}}   + \frac{1^{\frac{1}{3}}}{m_2 \sqrt[3]{m_1}} + \frac{1^{\frac{1}{3}}}{m_3 \sqrt[3]{m_1}}=0.     \]
The above equation holds if and and only if 
\[   m_1 \sqrt[3]{-(m_2+m_3)}  = m_2 \sqrt[3]{m_1}= m_3 \sqrt[3]{m_1}, \ {\it or}\ \frac{1}{m_1 \sqrt[3]{-(m_2+m_3)}} +\frac{1}{m_2 \sqrt[3]{m_1}}+\frac{1}{m_3 \sqrt[3]{m_1}}=0.   \]
The first case has no real solutions, and the second case reduces to 
\begin{equation} \label{equ:m-D4}\notag 
 \mathcal{V}_{V}[1,23]:~~~~~~~~~ m_1^2(m_2+m_3)^4=m_2^3m_3^3, 
\end{equation}
 which  has real solutions only if $m_2m_3>0$.

Note that 
\[  \sum_{j=1}^{4}m_j w_1 + m_4 w_{14} =   \sum_{j=1}^4 m_j w_j + m_2 w_{21} +m_3 w_{31}= m_2 w_{21} +m_2 w_{31} \preceq \epsilon,\] 
then $w_{14}\approx \epsilon^{-2}$ and $z_{14}\approx \epsilon^{2}$. 
It follows that 
\begin{equation}\label{DiagramV0}
 \sum_{j=1}^3 m_j\neq 0,  r_{12},r_{23},r_{13}\approx \epsilon, r_{14}\approx 1,  r_{24}\approx  r_{34}\approx \sqrt{\epsilon}^{-1}.
\end{equation}

Similarly,  Diagram V  has 11 other variants. Some of them are presented in the fifth row of Figure \ref{fig:Problematicdiagramsr21cj}. We still call them Diagram V if there is no confusion arise.  The corresponding mass sets are denoted by 
\[  \mathcal{V}_{V}[i, jk], \]
if  the fully edged triangle has 
vertices \textbf{i},  \textbf{j},  \textbf{k} and vertex \textbf{i} is connected to the fourth vertex by one $w$-stroke.  For instance,  for the   first two  diagrams in  the fifth row of Figure \ref{fig:Problematicdiagramsr21cj},  the corresponding mass sets are denoted by  $ \mathcal{V}_{V}[3,14]$ and  $ \mathcal{V}_{V}[4,13]$ respectively.  Both of the first two  diagrams will be referred  as Diagram V with $\triangle_{134}$.

\subsection{Diagram VI}
We will reach our result without discussing the mass polynomial of this diagram. Here is one remarks about this diagram.
Note that the momentum of inertia $I$ tends to zero by estimation 2. By Lemma \ref{lem:finitepotential}, the momentum of inertia is constant on a continuum of central configurations, so such a singular sequence exists only on the subset $I=0$. 

\begin{proposition}\label{prp:nonadjacent}
For singular sequences corresponding to Diagram V and VI, we have the following estimates
\[ r_{jk}r_{lp}\prec 1, \ r_{jk}^2r_{lp}\preceq 1,   \]
	where $r_{jk}$ and $r_{lp}$ are any two nonadjacent distances.  The same estimate holds for Diagram III, except for  third  case  of Subsection \ref{subsec:III-3}, for which it is necessary that  the masses satisfy equation \eqref{equ:m-D3b}. 
\end{proposition}

\section{Finiteness of Central configurations with $m\neq 0$ 
} \label{sec:pr1}

In this section, we esabilish the finiteness of central configurations in the situation where the toatl mass  $m\neq 0$. There are two cases, whether there is some partial sum of three masses being zero or not. 

The first case,  where all  partial sum of three masses are nonzero,  is considered in the first two subsections and we show that  

\begin{theorem}\label{thm:main-mn01}
	Suppose that $m_1, m_2, m_3, m_4$ are real  and nonzero,  if $\sum_{i=1}^4 m_i\ne0$ and $\prod_{1\leq j < k < l\leq 4}(m_j+m_k+m_l)\neq 0$, then   system \eqref{equ:sys-1}, which defines the normalized central configurations in the complex domain, possesses  finitely many solutions.
\end{theorem}

The second case,  where at least one  partial sum of three masses is  zero,  is considered 
in the next two subsections and we show that

 \begin{theorem}\label{thm:main-mnn01}
	Suppose that $m_1, m_2, m_3, m_4$ are real  and nonzero,  if $\sum_{i=1}^4 m_i\ne0$ and $\prod_{1\leq j < k < l\leq 4}(m_j+m_k+m_l)= 0$, then   system \eqref{equ:sys-1}, which defines the normalized central configurations in the complex domain, possesses  finitely many solutions.
\end{theorem}

Since the momentum of inertia is a constant on a continuum of solution of system \eqref{equ:sys-1},  the proof  of the above two  theorems completes after we establish Theorem \ref{thm:main-mn011},  Theorem \ref{thm:main-63}, Theorem \ref{thm:main-mnn011} and Theorem \ref{thm:main-mn012}.

\subsection{Finiteness with $m\neq 0$, $\prod_{1\leq  k\leq 4}(m-m_k)\neq 0$ and $I \neq 0$ }

 \begin{theorem}\label{thm:main-mn011}
 	Suppose that $m_1, m_2, m_3, m_4$ are real  and nonzero,  if $\sum_{i=1}^4 m_i\ne0$, $\prod_{1\leq j < k < l\leq 4}(m_j+m_k+m_l)\neq 0$ and $I \neq 0$, then   system \eqref{equ:sys-1}, which defines the normalized central configurations in the complex domain, possesses  finitely many solutions.
 \end{theorem}

 In this case only Diagram I, Diagram II, Diagram III and Diagram V are possible. From now on, we postpone the tedious computations of the proofs to the Appendix, and we will say $\triangle_{jkl}\longrightarrow0$, if we have one singular sequence corresponding to one diagram  with a fully edged triangle with  vertices  $\textbf{j, k, l}$.   \\

\noindent\emph{Proof of Theorem \ref{thm:main-mn011}}:
It is easy to see that giving five of $r_{kl}^2$'s, $1\le k<l\le 4$,  determines only finitely many geometrical  configurations up to rotation.
Suppose that there are infinitely many solutions of system (4)  in the complex domain.  Then  at least two of $r_{kl}^2$'s  must take infinitely many values and thus are dominating by Lemma \ref{lem:domin}.
Suppose that $r_{kl}^2$ is dominating for some  $1\le k<l\le 4$.  There must exist a  singular sequence of central configurations with $r_{kl}^2\longrightarrow0$,  which happens only in Diagram III or Diagram V. In either case, by  Proposition \ref{prp:nonadjacent},
\[   r_{kl}^2r_{ij}^2 \longrightarrow0, \  r_{ki}^2r_{jl}^2\longrightarrow0, \ r_{il}^2r_{kj}^2\longrightarrow0,  \]
along the singular sequence.  Then  all  the three polynomials
\[ r_{13}^2 r_{24}^2, \ r_{14}^2 r_{23}^2, \ r_{12}^2 r_{34}^2,  \]
are dominating.
It is easy to see that there exist $\triangle_{ijk}\longrightarrow0$ (say $\triangle_{123}\longrightarrow0$), so  $r^2_{12},r^2_{23}$ and $r^2_{13}$   are dominating.

There also  exist singular sequences with $r_{12}^2 r_{34}^2\longrightarrow\infty$, singular sequences with $r_{13}^2 r_{24}^2\longrightarrow\infty$ and singular sequences with $r_{14}^2 r_{23}^2\longrightarrow\infty$. These sequences must correspond to Diagram I or Diagram II.

Consider  a  singular sequence with $r_{13}^2 r_{24}^2\longrightarrow\infty$, thus, the masses must be belong to
\begin{equation}\label{equ:cond1}\notag
\mathcal{V}_{IA}[12,34]\bigcup  \mathcal{V}_{IA}[14,23]\bigcup  \mathcal{V}_{IB}[12,34]\bigcup  \mathcal{V}_{IB}[14,23]\bigcup  \mathcal{V}_{II}[13,24].
\end{equation}
Note that  $\mathcal{V}_{IA}[12,34]\bigcup \mathcal{V}_{IA}[14,23]\subset \mathcal{V}_{II}[13,24]$. Thus we the above set is the union of 
\begin{equation}\label{equ:newcond1}
\begin{array}{lc}
 1). ~~~\mathcal{V}_{II}[13,24];& \\
 2).~~~ \mathcal{V}_{IB}[12,34]\bigcup  \mathcal{V}_{IB}[14,23]. &
\end{array}
\end{equation}

Repeat the argument with  $r_{14}^2 r_{23}^2$. Then the  masses must belong to one of 
\begin{equation}\label{equ:newcond2}
\begin{array}{lc}
 1). ~~~\mathcal{V}_{II}[14,23];& \\
 2).~~~ \mathcal{V}_{IB}[12,34]\bigcup \mathcal{V}_{IB}[13,24]. &
\end{array}
\end{equation}

Repeat the argument  with  $r_{12}^2 r_{34}^2$. Then the   masses must belong to one of 
\begin{equation}\label{equ:newcond3}
\begin{array}{lc}
 1). ~~~\mathcal{V}_{II}[12,34];& \\
 2).~~~ \mathcal{V}_{IB}[13,24]\bigcup \mathcal{V}_{IB}[14,23]. &
\end{array}
\end{equation}

Our strategy is to  show that the three constraints \eqref{equ:newcond1}, \eqref{equ:newcond2} and \eqref{equ:newcond3}, together with other available constraints,  can not be satisfied simultaneously. 
 
 Define 
 \begin{align*}
 V[123]\triangleq & \mathcal{V}_{III}[123] \bigcup \mathcal{V}_{V}[1,23]\bigcup \mathcal{V}_{V}[2,13]\bigcup \mathcal{V}_{V}[3,12];\\
 V[124]\triangleq & \mathcal{V}_{III}[124] \bigcup \mathcal{V}_{V}[1,24]\bigcup \mathcal{V}_{V}[2,14]\bigcup \mathcal{V}_{V}[4,12];\\
 V[134]\triangleq & \mathcal{V}_{III}[134] \bigcup \mathcal{V}_{V}[1,34]\bigcup \mathcal{V}_{V}[3,14]\bigcup \mathcal{V}_{V}[4,13];\\
 V[234]\triangleq & \mathcal{V}_{III}[234] \bigcup \mathcal{V}_{V}[2,34]\bigcup \mathcal{V}_{V}[3,24]\bigcup \mathcal{V}_{V}[4,23]. 
 \end{align*}
{\bfseries{Case 1: Three of $1)$  are satisfied. }}
That is, $m_1m_3=m_2m_4, m_1m_4=m_2m_3, m_1m_2=m_3m_4$. Then $m_1^2=m_2^2=m_3^2=m_4^2$. Then  $m_1=m_2=m_3=m_4$. 
However, recall that $\triangle_{123}\longrightarrow0$,  consequently, the  masses belong to the set $V[123]$, a contradiction.

{\bfseries{Case 2: Two of $1)$  are satisfied}.}
Without lose of generality, assume that the masses  belong to the set $\mathcal{V}_{II}[13,24] \mathcal{V}_{II}[14,23]$, i.e.,  $m_1m_3=m_2m_4, m_1m_4=m_2m_3$.
Then  $m_1^2=m_2^2, m_3^2=m_4^2$. Then  $m_1=m_2, m_3=m_4$. 
The masses also belong to the set
\[   \mathcal{V}_{IB}[13,24]\bigcup \mathcal{V}_{IB}[14,23].\]
However, recall that $\triangle_{123}\longrightarrow0$, consequently, the  masses belong to the set $V[123]$, a contradiction. 

{\bfseries{Case 3: One of $1)$  is satisfied}.}
Without lose of generality, assume that the masses  belong to the set   $\mathcal{V}_{II}[13,24]$, i.e.,  $m_1m_3=m_2m_4$.
Then  the masses also belong to the set 
\[   \begin{array}{c}
      (\mathcal{V}_{IB}[12,34]\bigcup \mathcal{V}_{IB}[13,24]) (\mathcal{V}_{IB}[13,24]\bigcup \mathcal{V}_{IB}[14,23]) \\
       =(\mathcal{V}_{IB}[12,34]\mathcal{V}_{IB}[14,23])\bigcup\mathcal{V}_{IB}[13,24]. 
     \end{array}
\]

{\bfseries{Subcase 1: $\mathcal{V}_{IB}[12,34] \mathcal{V}_{IB}[14,23]$}}:

Note that $\triangle_{123}\longrightarrow0$. It is easy to check that $\mathcal{V}_{II}[13,24]\mathcal{V}_{IB}[12,34] \mathcal{V}_{IB}[14,23] V[123]$ is empty,  a contradiction.

\begin{figure}[h!] 
	\centering
	\includegraphics[width=7.2 cm, height=3 cm]{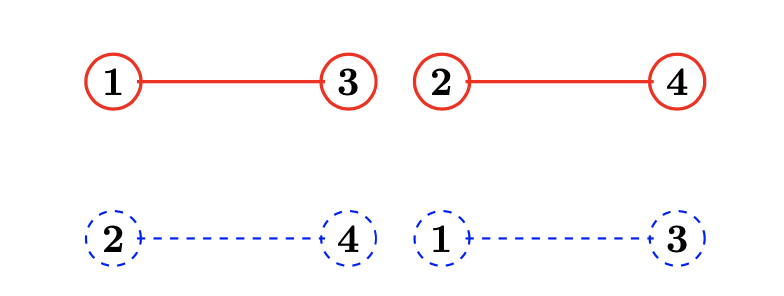}
	\caption{}
	\label{fig:sec6.1}
\end{figure}

{\bfseries{Subcase 2:  $ \mathcal{V}_{IB}[13,24]$}}:

Except that $\triangle_{123}\longrightarrow0$, note that in this case one of the two  diagrams in Figure \ref{fig:sec6.1}  occurs, thus $r_{14}^2$ and $r_{34}^2$ are dominating polynomials. Consider a singular sequence with $r^2_{14}\longrightarrow0$, then the sequence corresponds to Diagram III or V.  In Diagram III, the distances not in the fully edged triangle $\succeq 1$ since $I\ne 0$. Thus, we have 
$$\triangle_{124} \longrightarrow0, \ {\it or} \  \triangle_{134} \longrightarrow0. $$ 
Repeat the argument with $r_{34}^2$, we have 
$$\triangle_{134} \longrightarrow0, \ {\it or} \  \triangle_{234} \longrightarrow0. $$  
Thus, we conclude that we have 
$$\triangle_{134} \longrightarrow0, \ {\it or\ } \ \triangle_{124}\longrightarrow0,   \triangle_{234} \longrightarrow0. $$  
However, it is easy to check that both of the following two sets 
\begin{align*}
\mathcal{V}_{II}[13,24]\mathcal{V}_{IB}[13,24] V[123] V[134], \ 
\mathcal{V}_{II}[13,24]\mathcal{V}_{IB}[13,24] V[123] V[124] V[234], 
\end{align*}
are empty, a contradiction. 

{\bfseries{Case 4: No $1)$  is satisfied}.}Then  the masses belong to the set
\[  \begin{array}{c}
      (\mathcal{V}_{IB}[12,34]\bigcup  \mathcal{V}_{IB}[14,23]) (\mathcal{V}_{IB}[12,34]\bigcup \mathcal{V}_{IB}[13,24])(\mathcal{V}_{IB}[13,24]\bigcup \mathcal{V}_{IB}[14,23]) \\
      = (\mathcal{V}_{IB}[12,34] \mathcal{V}_{IB}[14,23])\bigcup (\mathcal{V}_{IB}[12,34]\mathcal{V}_{IB}[13,24])\bigcup(\mathcal{V}_{IB}[13,24] \mathcal{V}_{IB}[14,23]). 
    \end{array}
 \]
 
Without lose of generality, assume that the masses belong to the set $\mathcal{V}_{IB}[12,34] \mathcal{V}_{IB}[13,24]$. Except that $\triangle_{123}\longrightarrow0$, note that:

By $\mathcal{V}_{IB}[12,34]$, repeat the argument utilized above, we have two dominating polynomials,   $r_{14}^2,  r_{24}^2$.   Then let $r_{14}^2\longrightarrow 0$, we have $\triangle_{124} \longrightarrow0, \ {\it or} \  \triangle_{134} \longrightarrow0. $ Let $r_{24}^2\longrightarrow 0$, we have $\triangle_{124} \longrightarrow0, \ {\it or} \  \triangle_{234} \longrightarrow0. $  Thus, we conclude that we have 
$$\triangle_{124} \longrightarrow0, \ {\it or\ } \ \triangle_{134}\longrightarrow0,   \triangle_{234} \longrightarrow0. $$

Repeat the argument for  $\mathcal{V}_{IB}[13,24]$ and the two polynomials $r_{14}^2, r_{34}^2$.  We have 
$$\triangle_{134} \longrightarrow0, \ {\it or\ } \ \triangle_{124}\longrightarrow0,   \triangle_{234} \longrightarrow0. $$

In short, at least two of $\triangle_{124}\longrightarrow0$, $\triangle_{134}\longrightarrow0$ and $\triangle_{234}\longrightarrow0$ occur.

{\bfseries{Subcase 1: $\triangle_{124}\longrightarrow0 \ {\rm and}\ \triangle_{134}\longrightarrow0$}}:

However, it is easy to check that 
\[  \mathcal{V}_{IB}[12,34] \mathcal{V}_{IB}[13,24] V[123] V[124] V[134] \]
is empty, a contradiction. 

{\bfseries{Subcase 2: $\triangle_{124}\longrightarrow0 \ {\rm and}\ \triangle_{234}\longrightarrow0$}}:

However, it is easy to check that 
\[  \mathcal{V}_{IB}[12,34] \mathcal{V}_{IB}[13,24] V[123] V[124] V[234] \]
is empty, a contradiction. 

{\bfseries{Subcase 3: $\triangle_{134}\longrightarrow0 \ {\rm and}\ \triangle_{234}\longrightarrow0$}}:

However, it is easy to check that 
\[  \mathcal{V}_{IB}[12,34] \mathcal{V}_{IB}[13,24] V[123] V[134] V[234] \]
is empty, a contradiction.

To summarize, we proved that the  system \eqref{equ:sys-1} possesses  finitely many solutions, if $\sum_{i=1}^4 m_i\ne0$, $\prod_{1\leq j < k < l\leq 4}(m_j+m_k+m_l)\neq 0$ and $I \neq 0$.

$~~~~~~~~~~~~~~~~~~~~~~~~~~~~~~~~~~~~~~~~~~~~~~~~~~~~~~~~~~~~~~~~~~~~~~~~~~~~~~~~~~~~~~~~~~~~~~~~~~~~~~~~~~~~~~~~~~~~~~~~~~~~~~~~~~~~~~~~~~~~~~~~~~~~\Box$\\

\subsection{Finiteness with $m\neq 0$, $\prod_{1\leq  k\leq 4}(m-m_k)\neq 0$ and $I = 0$ }

 \begin{theorem}\label{thm:main-63}
 	Suppose that $m_1, m_2, m_3, m_4$ are real  and nonzero,  if $\sum_{i=1}^4 m_i\ne0$, $\prod_{1\leq j < k < l\leq 4}(m_j+m_k+m_l)\neq 0$ and $I = 0$, then   system \eqref{equ:sys-1}, which defines the normalized central configurations in the complex domain, possesses  finitely many solutions.
 \end{theorem}

 In this case only Diagram I, Diagram II, Diagram III, Diagram V and Diagram VI are possible.\\

\noindent\emph{Proof of Theorem \ref{thm:main-63}}:
It is easy to see that giving five of $r_{kl}^2$'s, $1\le k<l\le 4$,  determines only finitely many geometrical  configurations up to rotation.
Suppose that there are infinitely many solutions of system (4)  in the complex domain.  Then  at least two of $\{r_{kl}^2\}$'s  must take infinitely many values and thus are dominating by Lemma \ref{lem:domin}.
Suppose that $r_{kl}^2$ is dominating for some  $1\le k<l\le 4$.  There must exist a  singular sequence of central configurations with $r_{kl}^2\longrightarrow0$,  which happens only in Diagram III, Diagram V or Diagram VI. In either case, by  Proposition \ref{prp:nonadjacent},
\[   r_{kl}^2r_{ij}^2 \longrightarrow0, \  r_{ki}^2r_{jl}^2\longrightarrow0, \ r_{il}^2r_{kj}^2\longrightarrow0,  \]
along the singular sequence.  Then  all  the three polynomials
\[ r_{13}^2 r_{24}^2, \ r_{14}^2 r_{23}^2, \ r_{12}^2 r_{34}^2,  \]
are dominating.
It is easy to see that there exist $\triangle_{ijk}\longrightarrow0$ (say $\triangle_{123}\longrightarrow0$), so  $r^2_{12}, r^2_{23}$ and $r^2_{13}$   are dominating.

Then  there also  exist singular sequences with $r_{12}^2 r_{34}^2\longrightarrow\infty$, singular sequences with $r_{13}^2 r_{24}^2\longrightarrow\infty$ and singular sequences with $r_{14}^2 r_{23}^2\longrightarrow\infty$. These sequences must correspond to Diagram I or Diagram II.

Consider  a  singular sequence with $r_{13}^2 r_{24}^2\longrightarrow\infty$, thus, the masses must belong to the set
\begin{equation}\label{equ:cond4} \notag
\mathcal{V}_{I0}[12,34] \bigcup \mathcal{V}_{I0}[14,23]\bigcup \mathcal{V}_{II0}[13,24],
\end{equation}
which is the union of 
\begin{equation}\label{equ:newcond4}
\begin{array}{l}
 1). ~~~\mathcal{V}_{II0}[13,24]; \\
 2).~~~\mathcal{V}_{I0}[12,34] \bigcup \mathcal{V}_{I0}[14,23].
\end{array}
\end{equation}

Repeat the argument with  $r_{14}^2 r_{23}^2$. Then the  masses must belong to one of the sets
\begin{equation}\label{equ:newcond5}
\begin{array}{l}
 1). ~~~\mathcal{V}_{II0}[14,23]; \\
 2).~~~\mathcal{V}_{I0}[12,34] \bigcup \mathcal{V}_{I0}[13,24].
\end{array}
\end{equation}

Repeat the argument  with  $r_{12}^2 r_{34}^2$. Then the masses    belong to one of the sets
\begin{equation}\label{equ:newcond6}
\begin{array}{l}
 1). ~~~\mathcal{V}_{II0}[12,34]; \\
 2).~~~\mathcal{V}_{I0}[13,24] \bigcup \mathcal{V}_{I0}[14,23].
\end{array}
\end{equation}

Our strategy is to  show that the three constraints \eqref{equ:newcond4}, \eqref{equ:newcond5} and \eqref{equ:newcond6}  are incompatible.\\

{\bfseries{Case 1: Three of $1)$  are satisfied. }}
A straightforward computation shows that all the conditions  are impossible.
\\

{\bfseries{Case 2: Two of $1)$  are satisfied}.}
Without lose of generality, assume that the masses    belong to  the set
$$\mathcal{V}_{II0}[13,24]\mathcal{V}_{II0}[14,23](\mathcal{V}_{I0}[13,24] \bigcup \mathcal{V}_{I0}[14,23]).$$

However, one checks that the set is  empty. \\

{\bfseries{Case 3: One of $1)$  is satisfied}.}
Without lose of generality, assume that the masses    belong to  the set$$\begin{array}{c}
    \mathcal{V}_{II0}[13,24](\mathcal{V}_{I0}[12,34] \bigcup \mathcal{V}_{I0}[13,24])(\mathcal{V}_{I0}[13,24] \bigcup \mathcal{V}_{I0}[14,23]) \\
    = \mathcal{V}_{II0}[13,24]\mathcal{V}_{I0}[12,34] \mathcal{V}_{I0}[14,23]\bigcup \mathcal{V}_{II0}[13,24]\mathcal{V}_{I0}[13,24].
  \end{array}
$$

However, a straightforward computation shows that that the set is  empty.

{\bfseries{Case 4: No $1)$  is satisfied}.}
Then the masses    belong to  the set
$$\begin{array}{c}
    (\mathcal{V}_{I0}[12,34] \bigcup \mathcal{V}_{I0}[14,23])(\mathcal{V}_{I0}[12,34] \bigcup \mathcal{V}_{I0}[13,24])(\mathcal{V}_{I0}[13,24] \bigcup \mathcal{V}_{I0}[14,23]) \\
    = (\mathcal{V}_{I0}[12,34]  \mathcal{V}_{I0}[14,23])\bigcup(\mathcal{V}_{I0}[12,34]  \mathcal{V}_{I0}[13,24])\bigcup(\mathcal{V}_{I0}[13,24]  \mathcal{V}_{I0}[14,23]).
  \end{array}
$$
Without lose of generality, assume that the masses    belong to  the set $\mathcal{V}_{I0}[12,34]  \mathcal{V}_{I0}[13,24]$. A straightforward computation shows that the set is  empty.

To summarize, we proved that the  system \eqref{equ:sys-1} possesses  finitely many solutions, if $\sum_{i=1}^4 m_i\ne0$, $\prod_{1\leq j < k < l\leq 4}(m_j+m_k+m_l)\neq 0$ and $I = 0$.

$~~~~~~~~~~~~~~~~~~~~~~~~~~~~~~~~~~~~~~~~~~~~~~~~~~~~~~~~~~~~~~~~~~~~~~~~~~~~~~~~~~~~~~~~~~~~~~~~~~~~~~~~~~~~~~~~~~~~~~~~~~~~~~~~~~~~~~~~~~~~~~~~~~~~\Box$\\

\subsection{Finiteness of Central configurations with $m\neq 0$,   $\prod_{1\leq  k\leq 4}(m-m_k)= 0$ and $I\ne 0$ }

 In this case  Diagram I, Diagram II, Diagram III, Diagram IV , Diagram V and Diagram VI are all possible. 
 First, we establish a result without regard to the value of $I$.

 \begin{lemma}\label{r21cj}
If some product $r_{jk}^4r_{lm}^2$ is not dominating on the closed algebraic subset $\mathcal A$. Then system \eqref{equ:sys-1} possesses  finitely many solutions.
\end{lemma}
\noindent\emph{Proof of Lemma \ref{r21cj}}:

If system \eqref{equ:sys-1} possesses  infinitely many solutions, without lose of generality, assume that $r_{12}^4r_{34}^2$ is not dominating on the closed algebraic subset $\mathcal A$. Then some level set $r_{12}^4r_{34}^2\equiv const\neq 0$ also includes infinitely many solutions of system \eqref{equ:sys-1}.

\begin{figure}[t!]
	\centering
	 \includegraphics[width=14 cm, height=14 cm]{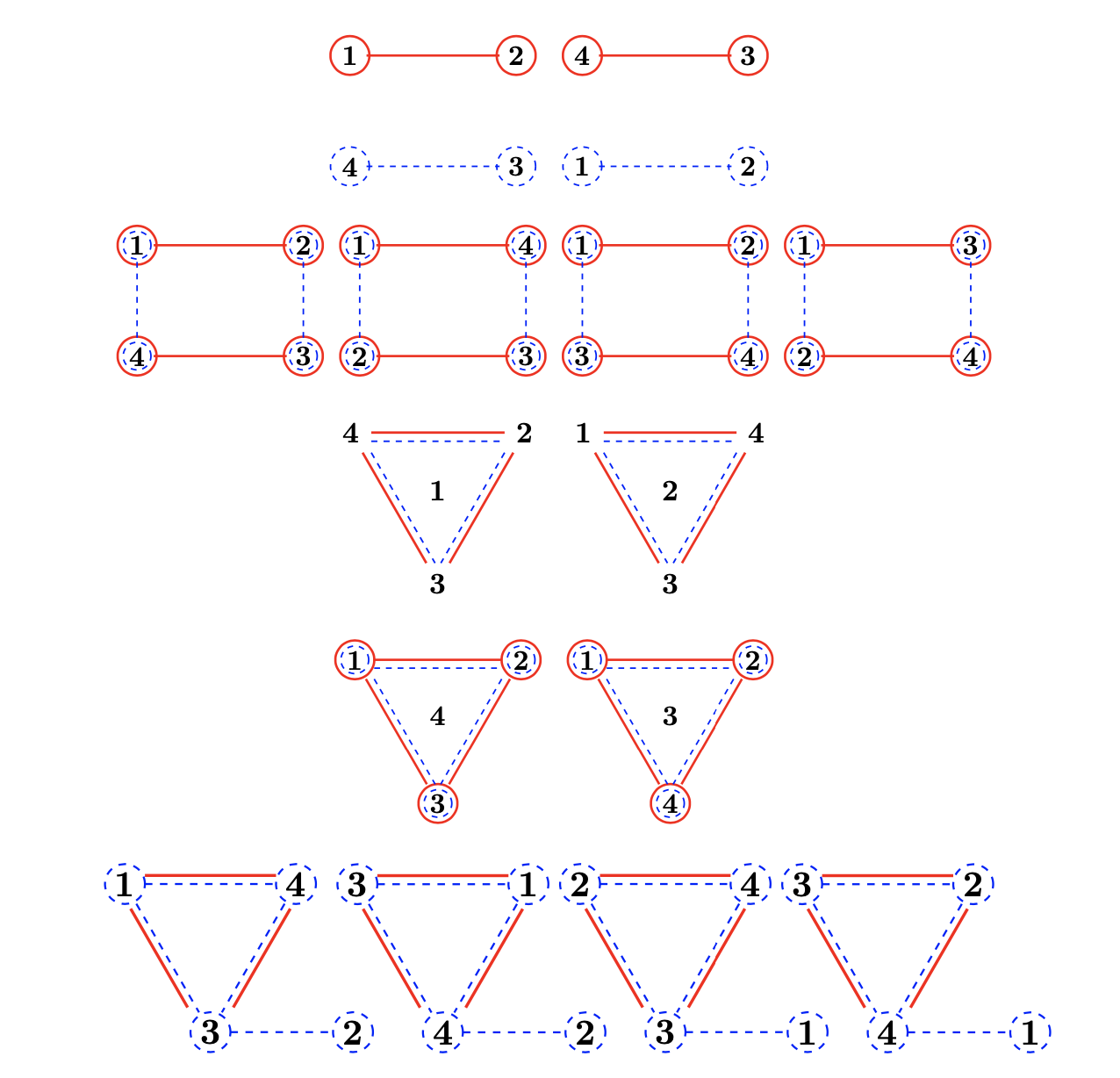}
\caption{All possible 
diagrams in the level set
 $r_{12}^4 r_{34}^2=const$}
\label{fig:Problematicdiagramsr21cj}
	
\end{figure}

In this level set, 
all  the corresponding 
possible diagrams  are presented in   Figure \ref{fig:Problematicdiagramsr21cj}. In particular, for Diagram III, we are in the third case of Subsection \ref{subsec:III-3}. More precisely,  if  the fully edged triangle  is $\triangle_{jkl}$  and the other vertex is   $\textbf{p}$, then 
$r^2_{pj}\sim r^2_{pk}\sim r^2_{pl}\approx \epsilon^{-1}$, and the masses belong to the set $\mathcal {V}_{III}[jkl]\mathcal {V}_{IV}[jkl]$.

Similarly,   at least two of $\{r_{kl}^2\}$'s  must take infinitely many values and thus are dominating.
Suppose that $r_{kl}^2$ is dominating for some  $1\le k<l\le 4$.  There must exist a  singular sequence of central configurations with $r_{kl}^2\longrightarrow0$,  which happens only in  Diagram III, Diagram IV or Diagram V. 
 By  Proposition \ref{prp:nonadjacent},  $r_{12}^2 r_{34}^2 \longrightarrow0$
along the singular sequence.  Then   the  polynomial $ r_{12}^2 r_{34}^2$
is  dominating.

Consider  a  singular sequence with $r_{12}^2 r_{34}^2\longrightarrow\infty$, then it is easy to see  that this happens only in  Diagram IV.  Without lose of generality, assume that Diagram IV with $\triangle_{123}$ occurs. Note that $r^4_{12} r^2_{14}r^2_{24}\longrightarrow\infty$ along a singular sequence corresponding to Diagram IV with $\triangle_{123}$, so $r^4_{12} r^2_{14}r^2_{24}$ is dominating. When 
$r^4_{12} r^2_{14}r^2_{24}\longrightarrow0$,  it is easy to see  that this happens only in Diagram IV with $\triangle_{124}$. Thus
\begin{equation}\label{r1223411}
    m_1+m_2+m_3=m_1+m_2+m_4=0, 
\end{equation}
and $r^2_{13}r^2_{14}$ is a dominating polynomial. 

Consider a singular sequence with $r^2_{13}r^2_{14}\longrightarrow0$. It is easy to see that this happens only in Diagram III with $\triangle_{134}$, and in the first two of Diagram V.  However,  it is easy to check that equation \eqref{r1223411} is not consistent with any one of the following three sets 
\[  \mathcal {V}_{III}[134]  \mathcal{V}_{IV}[134],\   \mathcal{V}_{V}[3,14], \  \mathcal{V}_{V}[4,13].  \]
This is a contradiction.
$~~~~~~~~~~~~~~~~~~~~~~~~~~~~~~~~~~~~~~~~~~~~~~~~~~~~~~~~~~~~~~~~~~~~~~~~~~~~~~~~~~~~~~~~~~~~~~~~~~~~~~~~~~~~~~~~~~~~~~~~~~~~~~~~~~~~~~~~~~~~~~~~~~~~\Box$\\

Hence, we assume that  all  $r_{jk}^4r_{lm}^2$'s  are dominating on the closed algebraic subset $\mathcal A$ from now on.

If  Diagram IV and Diagram III (the third case of Subsection \ref{subsec:III-3} with $r_{lp}^2 \succeq \epsilon^{-2}$, where vertex \textbf{p} is not in the fully edged triangle)  do not occur, then $r_{jk}r_{lp} \prec 1$ holds for  Diagram  III, 
 	and hence the proof of Theorem \ref{thm:main-mnn01} reduces to those of Theorem \ref{thm:main-mn01}.  Thus, we assume that either Diagram IV or Diagram III (the third case of Subsection \ref{subsec:III-3} with $r_{lp}^2 \succeq \epsilon^{-2}$, where vertex \textbf{p} is not in the fully edged  triangle) occurs. 
 Then there exists some $\triangle_{ijk}\longrightarrow 0$. 	
 	Without lose of generality, we assume 
 that $\triangle_{123}\longrightarrow0$ below. Then 
 \begin{equation}\label{IV123}
 m_1+m_2+m_3=0
 \end{equation}
 and it is easy to see that all of $r_{kl}^2$'s, $1\le k<l\le 4$, 
 are dominating.

Now we prove the finiteness in the case of $I\ne 0$. 

 \begin{theorem}\label{thm:main-mnn011}
 	Suppose that $m_1, m_2, m_3, m_4$ are real  and nonzero,  if $\sum_{i=1}^4 m_i\ne0$, $\prod_{1\leq j < k < l\leq 4}(m_j+m_k+m_l)= 0$ and $I \neq 0$, then   system \eqref{equ:sys-1}, which defines the normalized central configurations in the complex domain, possesses  finitely many solutions.
 \end{theorem}

 In this case only Diagram I, Diagram II, Diagram III, Diagram IV and Diagram V are possible. \\

\noindent\emph{Proof of Theorem \ref{thm:main-mnn011}}:

Consider a singular sequence with $r_{14}^2\longrightarrow 0$, 
 then the sequence corresponds to  $\triangle_{124}\longrightarrow0$, $\triangle_{134}\longrightarrow0$ in Diagram III, IV or V.  Repeat the argument with $r_{24}^2\longrightarrow 0$ and $r_{34}^2\longrightarrow 0$, then it is easy to see that
 at least two of $\triangle_{124}\longrightarrow0$, $\triangle_{134}\longrightarrow0$ and $\triangle_{234}\longrightarrow0$ occur.
Without lose of generality, assume that $\triangle_{124}\longrightarrow0$ and $\triangle_{134}\longrightarrow0$  occur.

Corresponding to $\triangle_{124}\longrightarrow0$ we have
\begin{equation}\label{Tri124}
    \begin{array}{c}
      (\frac{1}{ |m_1|^{\frac{1}{2}}} + \frac{1}{ |m_2|^{\frac{1}{2}}} +\frac{1}{ |m_4|^{\frac{1}{2}} } )( m_1+m_2+m_4)
      (m_1^2(m_2+m_4)^4-m_2^3m_4^3)\\
      (m_2^2(m_1+m_4)^4-m_1^3m_4^3)(m_4^2(m_2+m_1)^4-m_2^3m_1^3)=0
    \end{array}
\end{equation}
Corresponding to $\triangle_{134}\longrightarrow0$ we have
\begin{equation}\label{Tri134}
    \begin{array}{c}
      (\frac{1}{ |m_1|^{\frac{1}{2}}} + \frac{1}{ |m_3|^{\frac{1}{2}}} +\frac{1}{ |m_4|^{\frac{1}{2}} } )( m_1+m_3+m_4)
      (m_1^2(m_3+m_4)^4-m_3^3m_4^3)\\
      (m_3^2(m_1+m_4)^4-m_1^3m_4^3)(m_4^2(m_3+m_1)^4-m_3^3m_1^3)=0
    \end{array}
\end{equation}

Let us further consider $r_{12}^4r_{34}^2\longrightarrow \infty$,  then we are in Diagram I, Diagram II, Diagram III or Diagram IV.

{\bfseries{If we are in Diagram I of Diagram II. }}
Then the masses    belong to  the set
\begin{equation}\label{DiaIandII}
\mathcal{V}_{IA}[13,24] \bigcup\mathcal{V}_{IA}[14,23]\bigcup \mathcal{V}_{IB}[13,24] \bigcup\mathcal{V}_{IB}[14,23]\bigcup \mathcal{V}_{II}[12,34].
\end{equation}

However, a straightforward computation shows that  masses  of the above set are not  consistent with  equations (\ref{IV123}), (\ref{Tri124}) and (\ref{Tri134}).
\\

{\bfseries{If we are in Diagram III}.}
Then we have 
\begin{equation}\label{III134}
\left\{
             \begin{array}{lr}
             \frac{1}{ |m_1|^{\frac{1}{2}}} + \frac{1}{ |m_3|^{\frac{1}{2}}} +\frac{1}{ |m_4|^{\frac{1}{2}} }=0&  \\
             m_1+m_3+m_4=0 &
             \end{array}
\right.
\end{equation}
or
\begin{equation}\label{III234}
\left\{
             \begin{array}{lr}
             \frac{1}{ |m_2|^{\frac{1}{2}}} + \frac{1}{ |m_3|^{\frac{1}{2}}} +\frac{1}{ |m_4|^{\frac{1}{2}} }=0&  \\
             m_2+m_3+m_4=0 &
             \end{array}
\right.
\end{equation}

However, a straightforward computation shows that all of these relations are not  consistent with  equations (\ref{IV123}), (\ref{Tri124}) and (\ref{Tri134}).\\

{\bfseries{ If we are in Diagram IV}.}
Then it is Diagram IV with $\triangle_{234}$ or $\triangle_{134}$.

{\bfseries{Case 1: If we are in Diagram IV with $\triangle_{234}$}}:
Then 
\begin{equation}\label{IV234}
m_2+m_3+m_4=0.
\end{equation}
A straightforward computation shows that the above equation  is not  consistent with  equations (\ref{IV123}), (\ref{Tri124}) and (\ref{Tri134}).

{\bfseries{Case 2: If we are in Diagram IV with $\triangle_{134}$ }}:
Then 
\begin{equation}\label{IV134}
m_1+m_3+m_4=0.
\end{equation}
In this case, equation (\ref{IV134}) is   consistent with  equations (\ref{IV123}), (\ref{Tri124}) and (\ref{Tri134}) and it is necessary to include more polynomials.  Let us further consider $r_{13}^4r_{24}^2\longrightarrow \infty$, we are in Diagram I, Diagram II, Diagram III and Diagram IV. 

{\bfseries{Subcase 1: If in Diagram I or Diagram II.}} 
 Then 
 \begin{equation}\label{DiaIandIIagain}
\mathcal{V}_{IA}[12,34] \bigcup\mathcal{V}_{IA}[14,23]\bigcup \mathcal{V}_{IB}[12,34] \bigcup\mathcal{V}_{IB}[14,23]\bigcup\mathcal{V}_{II}[13,24]. 
\end{equation}

However, a straightforward computation shows that  masses  of the above set are not   consistent with  equations (\ref{IV123}), (\ref{Tri124}) and (\ref{IV134}).\\

{\bfseries{Subcase 2: If in Diagram III.}}
We have 
\begin{equation}\label{III124}
\left\{
\begin{array}{lr}
\frac{1}{ |m_1|^{\frac{1}{2}}} + \frac{1}{ |m_2|^{\frac{1}{2}}} +\frac{1}{ |m_4|^{\frac{1}{2}} }=0&  \\
m_1+m_2+m_4=0 &
\end{array}
\right.
\end{equation}
or equation \eqref{III234}. However, a straightforward computation shows that  neither  equation \eqref{III124} nor equation \eqref{III234} 
are  consistent with  equations (\ref{IV123}), (\ref{Tri124}) and (\ref{IV134}).\\

{\bfseries{Subcase 3: If in Diagram IV.}}
 We have $\triangle_{234}\longrightarrow 0$ or $\triangle_{124}\longrightarrow 0$. If $\triangle_{234}\longrightarrow 0$, then we have 
\[  m_2+m_3+m_4=0, \]
not consistent with  equations (\ref{IV123}), (\ref{Tri124}) and (\ref{IV134}).\\

If  $\triangle_{124}\longrightarrow 0$, then 
$ m_1+m_2+m_4=0.$
By \eqref{IV123} and \eqref{IV134}, the masses are 
\begin{equation}\label{equ:thm73-mass}
m_1=-2m_2, m_3= m_2, m_4= m_2. 
\end{equation}
Note that $r_{23}^2 r_{24}^2r_{34}^2\longrightarrow \infty$, it is dominating. Let $r_{23}^2 r_{24}^2r_{34}^2\longrightarrow 0$, we are in Diagram III, IV and V. The masses satisfying equation \eqref{equ:thm73-mass} can not admit Diagram III and V. In Diagram IV, the fully edged triangle can only be $\triangle_{234}$, so we have 
\[  m_2+m_3+m_4=0,\] 
which contradicts with \eqref{equ:thm73-mass}. 

To summarize, we proved that the  system \eqref{equ:sys-1} possesses  finitely many solutions, if $\sum_{i=1}^4 m_i\ne0$, $\prod_{1\leq j < k < l\leq 4}(m_j+m_k+m_l)= 0$ and $I \neq 0$.

$~~~~~~~~~~~~~~~~~~~~~~~~~~~~~~~~~~~~~~~~~~~~~~~~~~~~~~~~~~~~~~~~~~~~~~~~~~~~~~~~~~~~~~~~~~~~~~~~~~~~~~~~~~~~~~~~~~~~~~~~~~~~~~~~~~~~~~~~~~~~~~~~~~~~\Box$\\

\subsection{Finiteness with $m\neq 0$, $\prod_{1\leq  k\leq 4}(m-m_k)= 0$ and $I = 0$ }

 \begin{theorem}\label{thm:main-mn012}
 	Suppose that $m_1, m_2, m_3, m_4$ are real  and nonzero,  if $\sum_{i=1}^4 m_i\ne0$, $\prod_{1\leq j < k < l\leq 4}(m_j+m_k+m_l)\neq 0$ and $I = 0$, then   system \eqref{equ:sys-1}, which defines the normalized central configurations in the complex domain, possesses  finitely many solutions.
 \end{theorem}

 In this case all of Diagram I, Diagram II, Diagram III, Diagram IV, Diagram V and Diagram VI are possible.

\noindent\emph{Proof of Theorem \ref{thm:main-mn012}}:

Let us first consider $r_{12}^4r_{34}^2\longrightarrow \infty$,  then we are in Diagram I, Diagram II, Diagram III or Diagram IV.

Then the masses    belong to  the set
\begin{equation}\label{r{12}2r{34}}
   \begin{array}{c}
    \mathcal{V}_{I0}[13,24]\bigcup \mathcal{V}_{I0}[14,23]\bigcup\mathcal{V}_{II0}[12,34] \\
      \bigcup(\mathcal{V}_{III}[134]\mathcal{V}_{IV}[134])\bigcup(\mathcal{V}_{III}[234]\mathcal{V}_{IV}[234])\bigcup \mathcal{V}_{IV}[134]\bigcup \mathcal{V}_{IV}[234] \triangleq V[12,34].
   \end{array}  
\end{equation}
Similarly, when considering $r_{13}^4r_{24}^2\longrightarrow \infty$ the masses    belong to  the set
\begin{equation}\label{r{13}2r{24}}
   \begin{array}{c}
    \mathcal{V}_{I0}[12,34]\bigcup \mathcal{V}_{I0}[14,23]\bigcup\mathcal{V}_{II0}[13,24] \\
      \bigcup(\mathcal{V}_{III}[124]\mathcal{V}_{IV}[124])\bigcup(\mathcal{V}_{III}[234]\mathcal{V}_{IV}[234])\bigcup \mathcal{V}_{IV}[124]\bigcup \mathcal{V}_{IV}[234] \triangleq V[13,24];
   \end{array}
\end{equation}when considering $r_{23}^4r_{14}^2\longrightarrow \infty$ the masses    belong to  the set
\begin{equation}\label{r{23}2r{14}}
   \begin{array}{c}
    \mathcal{V}_{I0}[13,24]\bigcup \mathcal{V}_{I0}[12,34]\bigcup\mathcal{V}_{II0}[14,23] \\
      \bigcup(\mathcal{V}_{III}[124]\mathcal{V}_{IV}[124])\bigcup(\mathcal{V}_{III}[134]\mathcal{V}_{IV}[134])\bigcup \mathcal{V}_{IV}[124]\bigcup \mathcal{V}_{IV}[134] \triangleq V[23,14].
   \end{array}
\end{equation}

Or\begin{equation}
   \begin{array}{c}
   V[12,34]=(\mathcal{V}_{I0}[13,24])\bigcup (\mathcal{V}_{I0}[14,23])\bigcup(\mathcal{V}_{II0}[12,34])\bigcup \mathcal{V}_{IV}[134]\bigcup \mathcal{V}_{IV}[234];\\
    V[13,24]=(\mathcal{V}_{I0}[12,34])\bigcup (\mathcal{V}_{I0}[14,23])\bigcup(\mathcal{V}_{II0}[13,24])\bigcup \mathcal{V}_{IV}[124]\bigcup \mathcal{V}_{IV}[234] ; \\
     V[23,14]=(\mathcal{V}_{I0}[13,24])\bigcup (\mathcal{V}_{I0}[12,34])\bigcup(\mathcal{V}_{II0}[14,23])
      \bigcup \mathcal{V}_{IV}[124]\bigcup \mathcal{V}_{IV}[134].
   \end{array}
\end{equation}
Other than the above relations, recall that we have assumed that  the masses belong to $\mathcal{V}_{IV}[123]$.

{\bfseries{Case 1: If  Diagram I occurs. }}Without lose of generality, assume that the masses    belong to  the set $\mathcal{V}_{I0}[13,24]$, then we claim that
\begin{equation}\label{DiaI1occ}
(\mathcal{V}_{I0}[13,24])\bigcap V[13,24] \bigcap \mathcal{V}_{IV}[123]=\emptyset.
\end{equation}
A straightforward computation shows the claim.
\\

{\bfseries{Case 2: If  Diagram II occurs. }}Without lose of generality, assume that the masses    belong to  the set $\mathcal{V}_{II0}[12,34]$, then we claim that
\begin{equation}\label{DiaI1occ}
(\mathcal{V}_{II0}[12,34])\bigcap V[13,24] \bigcap \mathcal{V}_{IV}[123]\bigcap V[23, 14] =\emptyset.
\end{equation}
A straightforward computation shows the claim.
\\

Therefore,  Diagram I and Diagram II can not occur and  the masses belong to the following set 
\[ \mathcal{V}_{IV}[123]\left(    \mathcal{V}_{IV}[124]   \mathcal{V}_{IV}[134]\bigcup   \mathcal{V}_{IV}[124] \mathcal{V}_{IV}[234]\bigcup \mathcal{V}_{IV}[134]   \mathcal{V}_{IV}[234]\right), \]
which contains only three solutions
\begin{center}
	$m_2=m_3=m_4=-\frac{1}{2}m_1,  \ m_1=m_3=m_4=-\frac{1}{2}m_2, \ m _1=m_2=m_4=-\frac{1}{2}m_3.$
\end{center}

It is easy to see that the first  solution  is not consistent with the constraints corresponding to Diagram I, Diagram II, Diagram III or Diagram V, therefore, only Diagram IV and Diagram VI are possible now.
 We consider the function $r_{12}^2r_{23}^2r_{13}^2r_{14}^2r_{24}^2r_{34}^2$, then it is easy to see that it is dominating. When considering $r_{12}^2r_{23}^2r_{13}^2r_{14}^2r_{24}^2r_{34}^2\rightarrow 0$, we are in Diagram VI, thus the function $r_{12}^2r_{13}^2r_{14}^2$ is dominating. When considering $r_{12}^2r_{13}^2r_{14}^2\rightarrow \infty$, we are in Diagram IV with $\triangle_{234}\longrightarrow0$. As a result, we have $m_2+m_3+m_4=0$, which is a contradiction.

Similarly, the other two solutions do not admit infinite solutions neither.

To summarize, we proved that the  system \eqref{equ:sys-1} possesses  finitely many solutions, if $\sum_{j=1}^4 m_j\ne0$, $\prod_{1\leq j < k < l\leq 4}(m_j+m_k+m_l)= 0$ and $I = 0$.

$~~~~~~~~~~~~~~~~~~~~~~~~~~~~~~~~~~~~~~~~~~~~~~~~~~~~~~~~~~~~~~~~~~~~~~~~~~~~~~~~~~~~~~~~~~~~~~~~~~~~~~~~~~~~~~~~~~~~~~~~~~~~~~~~~~~~~~~~~~~~~~~~~~~~\Box$\\

\section{Finiteness of Central configurations with $m= 0$} \label{sec:pr3}

\subsection{Problematic diagrams with vanishing total mass}

We prove that the first and third diagram  in  Figure \ref{fig:C=3} and all diagrams  in Figure  \ref{fig:C=41} and Figure  \ref{fig:C=5} are impossible.

{\bfseries{Case 1: The First and third diagram  in  Figure \ref{fig:C=3}}.}
First, it is easy to see that 
\begin{equation}
  \begin{array}{c}
    m_1+m_2=0, \\
    m_3+m_4=0, \\
    w_1\sim w_2\sim w_3\sim w_4\sim a\epsilon^{-2}, \\
    z_{12}\approx z_{34}\approx w_{12}\approx w_{34}\approx \epsilon.
  \end{array}
\end{equation}
It follows that
\begin{equation}\notag
      w_{12}= m_3(W_{32}-W_{31})+m_4(W_{42}-W_{41}),
\end{equation}
we claim that $z_1\sim z_2\sim z_3\sim z_4$ for  the first three diagrams  in  Figure \ref{fig:C=3}. Otherwise, we have $z_{jk}\approx \epsilon^{-2}$ and then  $W_{jk}=\frac{1}{z_{jk}^{\frac{3}{2}}w_{jk}^{\frac{1}{2}}}\prec  \epsilon^{2}$ for $j=3,4,k=1,2$. This  contradicts with $w_{12}\approx \epsilon$. 

For  the  third diagram, we have $z_2, z_4 \prec z_1 ,z_2$. This is  a contradiction,  so the third  diagram  does  not exist.   

For  the first  diagram,  we have
\begin{equation}\notag 
      z_{1}\sim m_2 Z_{21}, \ 
      z_{2}\sim m_1 Z_{12}+m_3 Z_{32}. 
\end{equation}
This is  a contradiction,  so the first   diagram  does  not exist.

{\bfseries{Case 2: The Diagrams  in  Figure \ref{fig:C=41}}.}
First, it is easy to see that the second and third  diagram are impossible.

For  the first  diagram, note that $r_{14}, r_{24}, r_{34}\succ \epsilon$. 
Without loss of generality, assume that $w_{14}\succ \epsilon$, then
 we have $w_{24}\sim w_{34}\sim w_{14}\succ \epsilon$. By
 $$0=\sum_{j=1}^4 m_jw_{4}= \sum_{j=1}^3 m_jw_{j4},  $$
it follows that $ \sum_{j=1}^3 m_j=0$,  a contradiction.

{\bfseries{Case 3: The Diagram  in  Figure \ref{fig:C=5}}.}
First, it is easy to see that $w_{14}\succ \epsilon$, then
we have $w_{24}\sim w_{34}\sim w_{14}\succ \epsilon$. Repeat the argument in the previous case. We then arrive at  $ \sum_{j=1}^3 m_j=0$,  a contradiction.

\begin{figure}[t!]
	\centering
	 \includegraphics[width=13 cm, height= 11 cm]{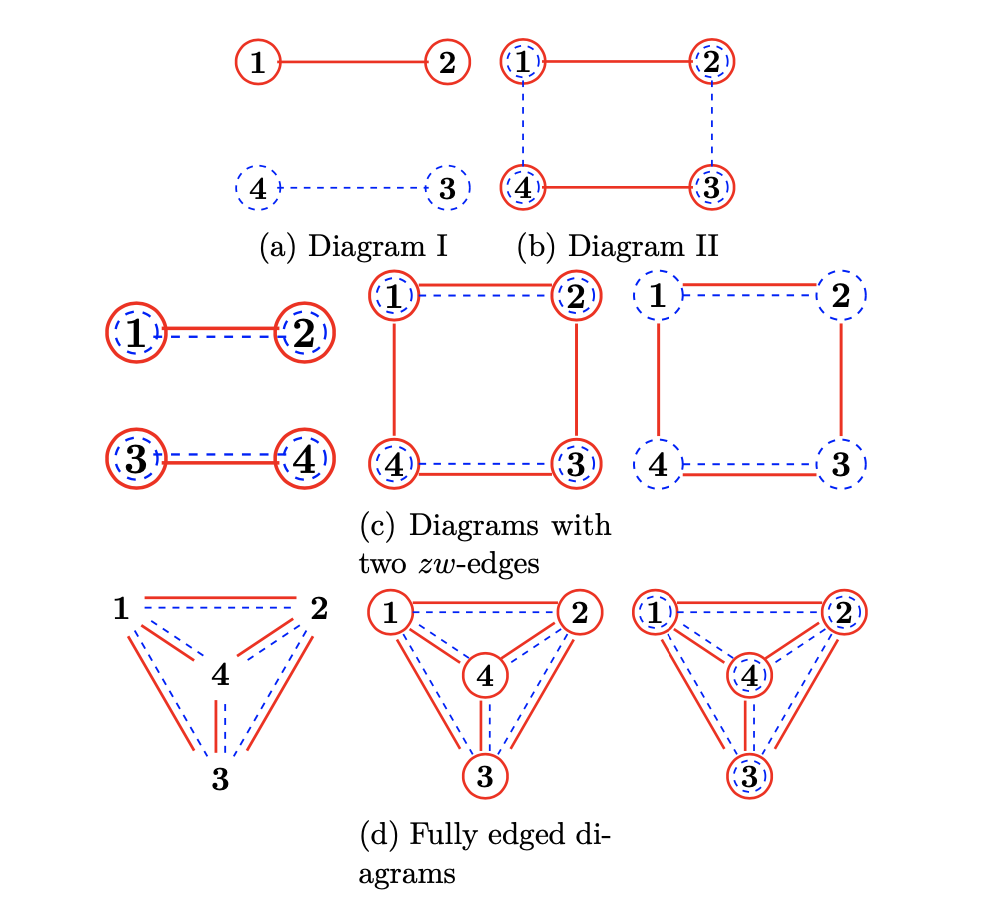}
	
	\caption{Problematic diagrams with $m=0$}
	\label{fig:Problematicdiagramsmis0}
\end{figure}

We could not eliminate the
diagrams in Figure \ref{fig:Problematicdiagramsmis0}. Some singular sequence could still exist and approach
any of these diagrams.

\subsubsection{Diagram I}

Note that the hypothesis that $m\ne 0$ is not utilized in Subsection \ref{subsec:DiagramI}, thus the estimates on the distances and mass polynomials derived there still holds for the present case. In particular, we have

\begin{equation}\label{D1m0mass} \notag
    (m_1+m_2) (m_3+m_4)\neq 0, 
\end{equation}

\begin{equation}\label{DiagramI03}
   r_{12}^3\sim m_1+m_2, \ 
      r_{34}^3\sim m_3+m_4, \    r_{jk}\approx \epsilon^{-2},
      j=1,2,k=3,4.
\end{equation}

The masses  corresponding to  Diagram I belong to one of the following two sets
\begin{center}
$\left\{
  \begin{array}{ll}
    \mathcal{V}_{IA0}[12,34]:, & \hbox{$m_1=m_2, \ m_3=m_4, \ m_1+m_3=0;$} \\
   \mathcal{V}_{IB}[12,34]:, & \hbox{$\mu_1 ^3 \mu_3^3  +  \mu_2 ^3 \mu_3^3 +  \mu_1 ^3 \mu_4^3 +  \mu_2 ^3 \mu_4^3=0,\  m_1m_2<0, \ m_3m_4<0,$}\\
     & \hbox{$ m_1+m_2\ne0, \ m_3+m_4\ne0.$}\\
  \end{array}
\right.$
\end{center}

\subsubsection{Diagram II}
Following the argument from \cite{Albouy2012Finiteness}, the masses  corresponding to Diagram II  satisfy the following equation
\begin{equation*} 
m_1m_3=m_2m_4.
\end{equation*}

Furthermore, we claim that
\begin{equation}\label{DiagramII02}
    m_1+m_2=m_2+m_3=m_3+m_4=m_1+m_4=0.
\end{equation}

Indeed, if $m_1+m_2\neq0$, then we have 
\[ z_{12}\approx \epsilon^{-2}, \ m_3+m_4\ne 0, \ z_{34}\approx \epsilon^{-2}.\] 
 By the fact $ z_{12}= (m_1+m_2) Z_{12}+m_3(Z_{32}-Z_{31})+m_4(Z_{42}-Z_{41})$, we have $z_{12}\sim  (m_1+m_2) Z_{12}$ and then $w_{12}\approx \epsilon^{2}$.  Hence,
$W_{12}=w_{12}^{-1/2}z_{12}^{-3/2}\approx \epsilon^2$.  Similarly, $W_{34}\approx \epsilon^2$.   We also would have $z_{13} \approx \epsilon^{-2}$ and $w_{13} \succ \epsilon$, then $W_{13} \prec \epsilon^2$. Similarly, $W_{24} \prec \epsilon^2$.

On the other hand, note $z_{13}\sim z_{12}, z_{14} \prec \epsilon.$
The the identity 
 $$0=\sum_{j=1}^4 m_jz_{1}= \sum_{j=2}^4 m_jz_{j1} $$
implies  that $ m_2+m_3=0$, then $ m_1+m_4=0$.

Hence, we have 
$w_{14}= (m_1+m_4) W_{14}+m_2(W_{24}-W_{21})+m_3(W_{34}-W_{31})  \preceq \epsilon^2$. By Estimate 1, $w_{14}\approx \epsilon^2$. Note that $W_{14}\approx \epsilon^{-2}$. It follows that $r_{14}\approx \epsilon^{4/3}$, which contradicts with Estimate 1. Thus equation \eqref{DiagramII02} follows.

Note that   \begin{equation}
  \begin{array}{c}
    r_{12},  r_{23}, r_{34}, r_{14} \prec 1, \\
     r_{13}\prec \epsilon^{-2},
    r_{24}\prec \epsilon^{-2}.
  \end{array}
\end{equation}

\subsubsection{Diagrams with two $zw$-edges}

Obviously, we have  $m_1+m_2 =m_3+m_4= 0$. Then the identities $0=\sum_{i=1}^4 m_i z_i =\sum_{i=1}^4 m_i w_i$ imply 
\begin{equation} \label{D2m0-1}
m_1 z_{12} =- m_3 z_{34}, \  m_1 w_{12} =- m_3 w_{34}. 
\end{equation}

For the last two diagrams with two $zw$-edges in Figure \ref{fig:Problematicdiagramsmis0}, we have 
\[   w_{23} =  (m_2+m_3) W_{23}+m_1(W_{13}-W_{12})+m_4(W_{43}-W_{42}) \prec \epsilon^{-2},  \]
which leads to 
\begin{equation} \label{D2m0-2}
m_1 W_{12} \sim m_3 W_{34}. 
\end{equation}
With \eqref{D2m0-1}, we have $m_1^2=m_3^2$.

For the first diagram with two $zw$-edges in Figure \ref{fig:Problematicdiagramsmis0}, we claim that at least one of $z_{23, }w_{23}$ shall have order less than $\epsilon^{-2}$. Otherwise, $W_{23}\approx W_{24}\approx W_{13}\approx W_{14} \approx \epsilon^4$. This contradicts with the equation 
\[ w_{12} =  (m_1+m_2) W_{12}+m_3(W_{32}-W_{31})+m_4(W_{42}-W_{41}).   \]
Then we have either 
\begin{equation} \label{D2m0-3}
m_1 W_{12} \sim m_3 W_{34}, \ {\it or} \ m_1 Z_{12} \sim m_3 Z_{34}. 
\end{equation}
With \eqref{D2m0-1}, we have $m_1^2=m_3^2$.

Therefore, the masses corresponding to Diagrams with two $zw$-edges satisfy   \begin{equation}\label{Diagramzw2}
m_1+m_2 =m_3+m_4= 0, \  m_1^2 =m_3^2. 
\end{equation}

\subsubsection{Fully edged diagrams}

We will reach our result without discussing the mass polynomial of this diagram. 

\subsection{Finiteness of Central configurations with $m= 0$}

 \begin{theorem}\label{thm:main-m01}
 	Suppose that $m_1, m_2, m_3, m_4$ are real  and nonzero with $\sum_{i=1}^4 m_j=0$, then   system \eqref{equ:sys-1}, which defines the normalized central configurations in the complex domain, possesses  finitely many solutions except perhaps  for, up to renumeration of the bodies, two groups of masses:  
 	$$m_1= m_2=- m_3=- m_4, \   m_1= m_2,  m_3=-\alpha^2 m_1,m_4=-(2-\alpha^2) m_1,$$
 	 where $\alpha$ is  the unique positive root of  equation  (\ref{xy}). 
 \end{theorem}

Numerically,   $\alpha \approx 1.2407$.

Note that if Diagram I with $\mathcal{V}_{IA0}$, Diagram II,  or Diagrams with two $zw$-edges occurs, then $m_k=\pm m_1, k=2,3,4$. 
Thus we assume that  Diagram I with $\mathcal{V}_{IA0}$, Diagram II, and  Diagrams with two $zw$-edges do not occur in the following discussion.\\

\noindent\emph{Proof of Theorem \ref{thm:main-m01}}:
It is easy to see that giving five of $r_{kl}^2$'s, $1\le k<l\le 4$,  determines only finitely many geometrical  configurations up to rotation.
Suppose that there are infinitely many solutions of system (4)  in the complex domain.  Then  at least two of $\{r_{kl}^2\}$'s  must take infinitely many values and thus are dominating by Lemma \ref{lem:domin}.
Suppose that $r_{kl}^2$ is dominating for some  $1\le k<l\le 4$.  There must exist a  singular sequence of central configurations with $r_{kl}^2\longrightarrow0$,  which happens only in the fully edged diagrams. In each case, it is easy to see that all of $r_{kl}^2$'s, $1\le k<l\le 4$, and $r_{13}^2 r_{24}^2, \ r_{14}^2 r_{23}^2, \ r_{12}^2 r_{34}^2$  are dominating.

Consider a singular sequence with  $r_{12}^2\longrightarrow \infty$, then the masses must belong to $\mathcal{V}_{IB}[13,24] \bigcup \mathcal{V}_{IB}[14,23]$. Consider a singular sequence with  $r_{13}^2\longrightarrow \infty$, then the masses must belong to $\mathcal{V}_{IB}[12,34] \bigcup \mathcal{V}_{IB}[14,23]$. Consider a singular sequence with  $r_{14}^2\longrightarrow \infty$, then the masses must belong to $\mathcal{V}_{IB}[12,34] \bigcup \mathcal{V}_{IB}[13,24]$. 
Hence,  at least two of $\mathcal{V}_{IB}[12,34]$, $\mathcal{V}_{IB}[13,24]$ and $\mathcal{V}_{IB}[14,23]$ occur.
Without lose of generality, assume that  $\mathcal{V}_{IB}[13,24]$ and $\mathcal{V}_{IB}[14,23]$  occur. That is, the masses belong to 
\[  \mathcal{V}_{IB}[13,24]\bigcap \mathcal{V}_{IB}[14,23] \bigcap \mathcal{V}_{m0},  \]
where $\mathcal{V}_{m0}$ denotes the set $\sum_{j=1}^4 m_j= 0$. 

It is easy to see  $m_1m_2>0$ and $m_3 m_4>0$.  Without loss of generality, assume that $m_1, m_2>0, m_3, m_4<0$.  Then a straightforward computation shows that the set $\mathcal{V}_{IB}[13,24]\bigcap \mathcal{V}_{IB}[14,23] \bigcap \mathcal{V}_{m0}$ is empty, provided that 
$ \prod_{i\ne j} (|m_i|-|m_j|)\ne 0.$ Note that  
\[ (m_1+m_3)(m_2+m_4)\ne0, \ (m_1+m_4)(m_2+m_3)\ne0,   \]
Hence, assume that $m_1=m_2$, or $m_3=m_4$.  A straightforward computation shows the set $\mathcal{V}_{IB}[13,24]\bigcap \mathcal{V}_{IB}[14,23] \bigcap \mathcal{V}_{m0}$ reduces to 
\begin{align*}
&\{ m_1=m_2, m_3 =-\alpha^2 m_1, m_4 =-(2-\alpha^2) m_1 \} \bigcup \{ m_1=m_2, m_3 =-(2-\alpha^2)  m_1, m_4 =-\alpha^2m_1 \}\\
&\bigcup  \{ m_1 =-\alpha^2 m_3, m_2 =-(2-\alpha^2) m_1, m_3=m_4 \} \bigcup \{ m_1=-(2-\alpha^2)  m_3, m_2 =-\alpha^2m_1,   m_3=m_4 \}
\end{align*}
where $\alpha$ is the unique positive root of the following equation:
 \begin{equation}\label{xy}
 (1-\alpha^3)^2-(1+\alpha^3)^2(2-\alpha^2)^3=0. 
 \end{equation}

$~~~~~~~~~~~~~~~~~~~~~~~~~~~~~~~~~~~~~~~~~~~~~~~~~~~~~~~~~~~~~~~~~~~~~~~~~~~~~~~~~~~~~~~~~~~~~~~~~~~~~~~~~~~~~~~~~~~~~~~~~~~~~~~~~~~~~~~~~~~~~~~~~~~~\Box$\\

\subsection{Finiteness with  $m_1=m_2=-m_3=-m_4$,  or  $m_1= m_2,  m_3=-\alpha^2 m_1,m_4=-(2-\alpha^2) m_1$. }

\begin{theorem}\label{thm:main-m02}
	Suppose that $m_1, m_2, m_3, m_4$ are real  and nonzero,  if $m_1=m_2=-m_3=-m_4$, or  $m_1= m_2,  m_3=-\alpha^2 m_1,m_4=-(2-\alpha^2) m_1$,  
	then   system \eqref{equ:sys-0} possesses  finitely many solutions.
\end{theorem}

We divide the proof into the collinear case and planar case. 

\subsubsection{The collinear case}

\noindent\emph{Proof of Theorem \ref{thm:main-m02}, Part 1}:

Assume that the masses are on the $x$-axis, and $x_i$ is the coordinate of $m_i$. Then 
system \eqref{equ:sys-0} becomes 
\begin{align}
 \sigma x_1 &= m_2 \frac{x_1-x_2}{|x_1-x_2|^{3} }+ m_3 \frac{x_1-x_3}{|x_1-x_3|^3 }+m_4\frac{x_1-x_4}{|x_1-x_4|^3 } \notag \\
&\cdots  \label{equ:sys-5}\\
 \sigma x_4 &= m_1 \frac{x_4-x_1}{|x_4-x_1|^3 } + m_2 \frac{x_4-x_2}{|x_4-x_2|^3 } +m_3\frac{x_4-x_3}{|x_4-x_3|^3 }, \notag 
\end{align}
where $\sigma =\pm 1$. 

If $m_1=m_2=-m_3=-m_4$, 
without lose of generality, assume $x_1<x_i$, $i=2, 3, 4$.   Since $\sum m_i x_i =0$, it suffices to assume that the ordering is  $x_1 < x_3 <x_4 <x_2$. A straightforward computation shows that  system  \eqref{equ:sys-5} has  finitely many  solutions.

 If $m_1= m_2,  m_3=-\alpha^2 m_1,m_4=-(2-\alpha^2) m_1$, by the reflection symmetry and the fact that $m_1=m_2$, it suffices to consider the following six orderings: $x_1<x_2<x_3<x_4$, $x_1<x_2<x_4<x_3$, $x_1<x_4<x_2<x_3$, $x_4<x_1<x_2<x_3$, $x_1<x_3<x_2<x_4$, $x_1<x_4<x_2<x_3$.  Note that  the first two orderings are impossible since $ \frac{ x_1 + x_2 }{2} \in (x_1, x_2)$,  and  $\frac {-m_3 x_3 - m_4 x_4 }{2}>x_2$. The third ordering is also impossible since
 \[   x_4 = \frac{x_1 + x_2+m_3 x_3}{-m_4}=  \frac{(2+m_3)x_1 + (1+m_3)(x_2-x_1)-m_3(x_2-x_3) }{-m_4}<x_1.   \]

 For each of the last three orderings, a straightforward computation shows that  system  \eqref{equ:sys-5} has  finitely many  solutions.  For instance, let $x_4<x_1<x_2<x_3$ and $\sigma=1$. The system becomes
\begin{align}
 x_1 + \frac{1}{d_{12}^2 }+  \frac{m_3}{d_{13}^2 }-\frac{(-2-m_3)}{d_{14}^2 }&=0\notag \\
x_2 - \frac{1}{d_{12}^2 }+  \frac{m_3}{d_{23}^2 }-\frac{(-2-m_3)}{d_{24}^2 }&=0\label{equ:sys-6} \\
x_3 - \frac{1}{d_{13}^2 }-  \frac{m_2}{d_{23}^2 }-\frac{(-2-m_3)}{d_{34}^2 }&=0, \notag
\end{align}
where $x_2=x_1+d_{12}, x_3=x_1+d_{13},x_4=x_1-d_{14}$, $d_{12}=-m_3 d_{13} -(2+m_3)d_{14}, d_{23}=d_{13}-d_{12}, d_{24}=d_{12}+d_{14}, d_{34}=d_{13}+d_{14}$.   After clearing the denominators, the above system consists of three polynomials in $x_1, d_{13}, d_{14}$.  Denote them by $p_1, p_2, p_3$ respectively. Note that each of them is  degree one in $x_1$. We  compute the resultants $r_1(d_{13}, d_{14})=Res(p_1, p_2, x_1)$ and $r_2(d_{13}, d_{14})=Res(p_1, p_3, x_1)$, and then compute $r_3(d_{14})=Res(r_1, r_3, d_{13})$, which is a nonzero polynomial in $d_{14}$ with degree 135. Thus, we conclude that   system  \eqref{equ:sys-6} has  finitely many  solutions.

$~~~~~~~~~~~~~~~~~~~~~~~~~~~~~~~~~~~~~~~~~~~~~~~~~~~~~~~~~~~~~~~~~~~~~~~~~~~~~~~~~~~~~~~~~~~~~~~~~~~~~~~~~~~~~~~~~~~~~~~~~~~~~~~~~~~~~~~~~~~~~~~~~~~~\Box$\\

\subsubsection{The planar case}

Recall that solving system \eqref{equ:sys-0}  is equivalent to finding critical points of the potential restricted on the set $I=c$.   For these  two groups of masses, following Celli, \cite{theis-Celli}, we use $r_{jk}$ as coordinates. For system with zero total mass,  there are extra constraints,  namely
\[  \sum_{k\ne 1} m_k r^2_{1k}=\sum_{k\ne 2} m_kr^2_{2k} =\sum_{k\ne 3} m_kr^2_{3k} =\sum_{k\ne 4} m_4r^2_{4k} =c_0.   \]
Then 
system \eqref{equ:sys-0} is equivalent to the following system
\begin{equation}\label{equivalentform1}
\begin{array}{cc}
m_jm_k+2(\lambda_j m_k+\lambda_k m_j)r_{jk}^3=0, & 1\leq j <k\leq 4 \\
\sum_{j\neq k}m_j r_{jk}^2=I, &  1\leq k\leq 4,
\end{array}
\end{equation}
where $\lambda_j$, $1\le j\le 4$, are the Lagrange multiplier.

\noindent\emph{Proof of Theorem \ref{thm:main-m02}, Part 2}:

Consider system \eqref{equivalentform1}.

{\bfseries{Case 1}}: $m_1=m_2=-m_3=-m_4. $

In this case, it is easy to see that
\[r_{12}=r_{14}=r_{23}=r_{24}\]
By  imposing the normalization $r_{24}=1$ and setting $ r_{34}=t$, system (\ref{equivalentform1}) reduces to 
\[   r_{12}^2=4-t^2, \ \frac{1}{(4-t^2)^{3/2}} +\frac{1}{t^3}=2. \]
Thus, $t\approx 0.816713, 1.825645$ and  the system has only two solutions.

{\bfseries{Case 2: $m_1= m_2,  m_3=-\alpha^2 m_1,m_4=-(2-\alpha^2) m_1$. }}
In this case, after eliminating the $\lambda_j$ and $I$
we have the relations:
\begin{equation}\label{equivalentform2} \notag 
\begin{array}{cc}
  \frac{1}{r_{12}^3}+\frac{1}{{r_{34}}^3}=\frac{1}{r_{13}^3}+\frac{1}{{r_{24}}^3}=\frac{1}{r_{14}^3}+\frac{1}{{r_{23}}^3}, &  \\
  \alpha^2 \left(r_{13}^2-r_{23}^2\right)=(2-\alpha^2) \left(r_{24}^2-r_{14}^2\right). & 
\end{array}
\end{equation}
It follows that 
\begin{equation}\notag 
    r_{13}-r_{23}=r_{14}-r_{24}=0,
\end{equation}
then by  imposing the normalization $r_{34}=1$, system  (\ref{equivalentform1}) reduces to 
\begin{equation}\label{equivalentform3}
\begin{array}{c}
   \frac{1}{r_{12}^3}=\frac{1}{r_{13}^3}+\frac{1}{r_{14}^3}-1,  \\
  r_{12}^2-\left(r_{13}^2-1\right) \alpha^2-r_{14}^2 \left(4-\alpha^2\right)=0,\\
  r_{13}^2- r_{14}^2+\alpha^2-1=0. 
\end{array}
\end{equation}
The last two equations of system \eqref{equivalentform3} imply $r_{12}^2= 4r_{14}^2 - \alpha^4$. Then $r_{14}$ satisfies 
$$ \frac{1}{(4 r_{14}^2 -\alpha^4)^{3/2}}=\frac{1}{(r_{14}^2 +1-\alpha^2)^{3/2}}+\frac{1}{r_{14}^3}-1. $$
which has finitely many  solutions.  
 Thus, system (\ref{equivalentform3}) has finitely many  solutions.

$~~~~~~~~~~~~~~~~~~~~~~~~~~~~~~~~~~~~~~~~~~~~~~~~~~~~~~~~~~~~~~~~~~~~~~~~~~~~~~~~~~~~~~~~~~~~~~~~~~~~~~~~~~~~~~~~~~~~~~~~~~~~~~~~~~~~~~~~~~~~~~~~~~~~\Box$\\

\newpage 
\section {Appendix}

Some  tedious computations in the proofs of Theorems in Section \ref{sec:pr1}  and  Section \ref{sec:pr3}  are omitted there. We now provide them in this appendix. 

 We denote  $\pm 1$ by $\sigma_i$, and $\sqrt{|m_i|}>0$  by  $\mu_i$, $i=1, 2, 3, 4$.   Without lose of generality, let $\mu_1=1$. Then those algebraic  constraints can be written as   polynomials  of the three real variable $\mu_2, \mu_3, \mu_4$. For instance, the constraints in $\mathcal{V}_{IB}[12, 34]$ are 
\begin{align*}
 &(\mu_3^3  + \sigma_1  \mu_4^3 )+ \sigma_2 \mu_2 ^3 (\mu_3^3  + \sigma_3  \mu_4^3 )=0, \\
 &   m_2<0, \ m_3m_4<0, \ 1+m_2\ne 0, \ m_3+m_4\ne 0. 
\end{align*}
We then see that $\sigma_1 \sigma_3=-1$ and the polynomial reduces to $$(\mu_3^3  + \sigma_1  \mu_4^3 )+ \sigma_2 \mu_2 ^3 (\mu_3^3  - \sigma_1  \mu_4^3 )=0.$$
 Note that $\sigma_1=1$ if and only if $1<\mu_2$.  The constraint in $\mathcal{V}_{III}[124]$ is 
 \[ ( \mu_2-\frac{\mu_4}{\mu_4-1})     ( \mu_2-\frac{\mu_4}{1-\mu_4}) ( \mu_2-\frac{\mu_4}{1+\mu_4}) =0.  \]

\emph{Computations in the proof of Theorem \ref{thm:main-mn011}}:

\textbf{Case 2}: The masses satisfy 
$m_1=m_2=1, m_3=m_4$, and it belongs to 
 \[ \mathcal{V}_{IB}[13,24] \bigcup \mathcal{V}_{IB}[14,23]  \bigcap  \left(\mathcal{V}_{III}[123] \bigcup ( \mathcal{V}_{V}[1,23]  \bigcup \mathcal{V}_{V}[2,13]\bigcup  \mathcal{V}_{V}[3,12] )\right). \] 
Then we have $m_1m_3<0$, so it can not belong to $\mathcal{V}_{V}[1,23]  \bigcup \mathcal{V}_{V}[2,13]$.  

Thus, the masses are  in $\mathcal{V}_{III}[123]\bigcup  \mathcal{V}_{V}[3,12]$, then $m_3=m_4=-\frac{1}{4}$. 
Obviously, the masses is not in  $\mathcal{V}_{IB}[13,24] \bigcup \mathcal{V}_{IB}[14,23]$. 

\textbf{Case 3, subcase 1}: It is easy to see that the masses are in 
\[  \mathcal{V}_{II}[13,24]\mathcal{V}_{IB}[12,34]  \mathcal{V}_{IB}[14,23] (\mathcal{V}_{III}[123] \bigcup \mathcal{V}_{V}[2,13] ).  \]
Then $\mu_3=\mu_2\mu_4$. By the fact that the masses are in $\mathcal{V}_{IB}[12,34]$, we find 
\[  \left(\mu_2^6+1\right)^2 \left(\mu_2^{12}-6 \mu_2^6+1\right) \mu_4^{12} =0. \]
Since $m_1+m_2\ne 0$, we obtain $\mu_2=(\sqrt{2}\pm1)^{1/3}$. Similarly, we have $\mu_4=(\sqrt{2}\pm1)^{1/3}$. It is easy to check that the masses are not in $\mathcal{V}_{III}[123] \bigcup \mathcal{V}_{V}[2,13]$.

\textbf{Case 3, subcase 2}: It is easy to see that the masses are in 
\[  \mathcal{V}_{II}[13,24]\mathcal{V}_{IB}[13,24]  (\mathcal{V}_{III}[123] \bigcup \mathcal{V}_{V}[1,23] \bigcup \mathcal{V}_{V}[3,12]) (\mathcal{V}_{III}[134] \bigcup \mathcal{V}_{V}[1,34] \bigcup \mathcal{V}_{V}[3,14]).  \]
By the fact that the masses are in $\mathcal{V}_{II}[13,24]$, we have $\mu_3=\mu_2\mu_4$. Substituting $\mu_3=\mu_2\mu_4$ into the polynomials corresponding to $\mathcal{V}_{IB}[13,24]$,  $\mathcal{V}_{III}[123] \bigcup \mathcal{V}_{V}[1,23] \bigcup \mathcal{V}_{V}[3,12]$, and $\mathcal{V}_{III}[134] \bigcup \mathcal{V}_{V}[1,34] \bigcup \mathcal{V}_{V}[3,14])$, we obtain three polynomials in terms of $(\mu_2, \mu_4)$, with degree 36, 17, and 17 respectively.  We will not write them explicitly here. Straightforward computation, for instance, by using Gr\"{o}bner basis, shows that they have no common positive root. Hence,  the set is empty. 

\textbf{Case 3, subcase 3}: It is easy to see that the masses are in 
\begin{align*}
& \mathcal{V}_{II}[13,24]\mathcal{V}_{IB}[13,24] \bigcap  (\mathcal{V}_{III}[123] \bigcup \mathcal{V}_{V}[1,23] \bigcup \mathcal{V}_{V}[3,12]) \bigcap \\
&(\mathcal{V}_{III}[124] \bigcup \mathcal{V}_{V}[2,14] \bigcup \mathcal{V}_{V}[4,12])\bigcap (\mathcal{V}_{III}[234] \bigcup \mathcal{V}_{V}[2,34] \bigcup \mathcal{V}_{V}[4,23]). 
\end{align*}
This set is also empty. The proof is similar to the above subcase and is omitted. 

\textbf{Case 4, subcase 1}: It is easy to see that the masses are in 
\begin{align*}
& \mathcal{V}_{IB}[12,34]\mathcal{V}_{IB}[13,24] \bigcap  (\mathcal{V}_{III}[123] \bigcup \mathcal{V}_{V}[1,23]) \bigcap \\
&(\mathcal{V}_{III}[124] \bigcup \mathcal{V}_{V}[2,14] )\bigcap (\mathcal{V}_{III}[134] \bigcup \mathcal{V}_{V}[3,14] ). 
\end{align*}

Note that the relation corresponding to $\mathcal{V}_{III}[124] \bigcup \mathcal{V}_{V}[2,14]$ can be written as
\[ 0=( \mu_2-\frac{\mu_4}{\mu_4-1})     ( \mu_2-\frac{\mu_4}{1-\mu_4}) ( \mu_2-\frac{\mu_4}{1+\mu_4}) ( \mu_2-\frac{\mu_4^{3/2}}{1+\mu_4^2}).      \]
There are four possibilities. 

If $\mu_2=\frac{\mu_4}{\mu_4-1}$, then $\mu_2>1$ and the constraints from $\mathcal{V}_{IB}[12,34]$ is
\[  \mu_2^6(\mu_3^3-\mu_4^3)^2-(\mu_3^3+\mu_4^3)^2=0. \] 
Substituting $\mu_2=\frac{\mu_4}{\mu_4-1}$ into it and clearing the denominator, we obtain one polynomial in terms $(\mu_3, \mu_4)$ with degree 12. Similarly, with the constraints from $\mathcal{V}_{III}[123] \bigcup \mathcal{V}_{V}[1,23]$ and  $\mathcal{V}_{III}[134] \bigcup \mathcal{V}_{V}[3,14]$, we obtain two polynomials in terms $(\mu_3, \mu_4)$ with degree 12.  Straightforward computation shows that they have no common positive root.  Hence, the set is empty if $\mu_2=\frac{\mu_4}{\mu_4-1}$. 

Similarly, we can show that the set is empty if $\mu_2=\frac{\mu_4}{1-\mu_4}$,  $\mu_2=\frac{\mu_4}{\mu_4+1}$, or $\mu_2=\frac{\mu_4^{3/2}}{1+\mu_4^2}$.

To summarize, the set corresponding to this subcase is empty.

\textbf{Case 4, subcase 2 and subcase 3}: It is easy to see that the masses are in 
\begin{align*}
& \mathcal{V}_{IB}[12,34]\mathcal{V}_{IB}[13,24] \bigcap  (\mathcal{V}_{III}[123] \bigcup \mathcal{V}_{V}[1,23]) \bigcap \\
&(\mathcal{V}_{III}[234] \bigcup \mathcal{V}_{V}[4,23] )\bigcap (\mathcal{V}_{III}[124] \bigcup \mathcal{V}_{V}[2,14] ), 
\end{align*}
and 
\begin{align*}
& \mathcal{V}_{IB}[12,34]\mathcal{V}_{IB}[13,24] \bigcap  (\mathcal{V}_{III}[123] \bigcup \mathcal{V}_{V}[1,23]) \bigcap \\
&(\mathcal{V}_{III}[234] \bigcup \mathcal{V}_{V}[4,23] )\bigcap (\mathcal{V}_{III}[134] \bigcup \mathcal{V}_{V}[3,14] ), 
\end{align*}
respectively. The two  sets are also empty. The proof is similar to the above subcase and is omitted. \\

\emph{Computations in the proof of Theorem \ref{thm:main-63}}：

\textbf{Case 1}: The masses are in $\mathcal{V}_{II0}[13,24] \mathcal{V}_{II0}[14,23]\mathcal{V}_{II0}[12,34]$. We have $m_2^2=m_3^2=m_4^2=1=m_1$ and $m_i+m_j\ne 0$. Then $m_1=m_2=m_3=m_4=1$, which does not belong to $\mathcal{I}_{II}[13,24]$, a contradiction. \\

\textbf{Case 2}: Without lose of generality, consider only the case $\mathcal{V}_{II0}[13,24]\mathcal{V}_{II0}[14,23]\mathcal{V}_{I0}[13,24]$. 
First we have 
\[ m_1=m_2^2=1, \ m_3^2=m_4^2, \ (m_1+m_2)(m_2+m_3)(m_3+m_4)\ne 0, \ m_1m_3<0, \ m_2m_4<0.  \]
Thus, we set the masses are $m_1=m_2=1, m_3=m_4=s<0, s\ne -1$, which does not belong to $\mathcal{I}_{I}[13,24]$, a contradiction. \\

\textbf{Case 3}: The masses are in either $ \mathcal{V}_{II0}[13,24]\mathcal{V}_{I0}[12,34] \mathcal{V}_{I0}[14,23]$ or $ \mathcal{V}_{II0}[13,24]\mathcal{V}_{I0}[13,24]$.

In the first subcase,  the signs of the masses are $(+, -, +, -)$. Set $m_2=s<0, m_4=t<0$. Then $m_3=st$.  Then  that the masses are in $\mathcal{I}_{I}[12,34]$ implies 
\[   \frac{s}{\sqrt[3]{1+s}} + \frac{st^2}{\sqrt[3]{t}\sqrt[3]{1+s} }= 0.  \]
Then $t=-1$ and $m_1+m_4=0$, hence  the masses are not in $\mathcal{V}_{IB}[14,23]$, a contradiction.

In the second subcase, we have  $m_1m_3=m_2m_4$. 
Then the masses are in $\mathcal{I}_{I}[13,24]$ implies 
\[  m_1+m_3 = -(m_2+m_4),  \]
 which contradicts with the hypothesis that the total mass $m\ne 0$.

\textbf{Case 4}:  The masses are in  $\mathcal{V}_{I0}[12,34]\mathcal{V}_{I0}[13,24]$.  
The masses are $(1, -\mu_2^2, -\mu_3^2, \mu_4^2)$. Then the masses satisfy the system 
\begin{align*}
&f_1=(\mu_3^3  + \sigma_1  \mu_4^3 )+ \sigma_2 \mu_2 ^3 (\mu_3^3  - \sigma_1  \mu_4^3 )=0, 
&f_2=(\mu_2^3  + \sigma_3  \mu_4^3 )+ \sigma_4 \mu_3 ^3 (\mu_2^3  - \sigma_3  \mu_4^3 )=0, \\
&\mu_2^6(\mu_3^2-\mu_4^2)-\mu_3^6\mu_4^6(1-\mu_2^2)=0, 
 &\mu_3^6(\mu_2^2-\mu_4^2)-\mu_2^6\mu_4^6(1-\mu_3^2)=0. 
\end{align*} 
The last two equations imply 
\[  \mu_3>\mu_4, \Leftrightarrow \mu_2<1, \ {\rm and \ }  \mu_2>\mu_4, \Leftrightarrow \mu_3<1. \]
Hence, there are four possibilities: $ \mu_2<1, \mu_3<1$; $ \mu_2<1, \mu_3>1$; $\mu_2>1, \mu_3>1$; $\mu_2>1, \mu_3<1$. Without lose of generality, we consider only the first two subcases. 

If $\mu_2<1, \mu_3<1$, then we have $\sigma_1=\sigma_2=\sigma_3=\sigma_4=-1$, i.e., 
\[  f_1=(\mu_3^3  - \mu_4^3 )- \mu_2 ^3 (\mu_3^3 +  \mu_4^3 )=0, \
f_2=(\mu_2^3  -   \mu_4^3 )-\mu_3 ^3 (\mu_2^3  + \mu_4^3 )=0. \]
Then we have 
\[ 0=f_{1}-f_{2}= (1+\mu_4^3) (\mu_3^3-\mu_2^3),  \]
thus $\mu_2=\mu_3$.  By the equation the third equation, we get $\mu_3=\mu_4\sqrt{\frac{1+\mu_4^4}{1+\mu_4^6}}$. Then 
\[  0=f_{1}= \mu_4^3 \left(\frac{\left(\mu_4^5+\mu_4\right)^3}{\left(\mu_4^6+1\right)^3}+\left(\mu_4^3-1\right) \left(\frac{\mu_4^4+1}{\mu_4^6+1}\right)^{3/2}+1\right). \]
It is easy to verify that the above equation has no  root for  $\mu_4\in(0, 1)$. Thus, the system has no positive solutions.

If $\mu_2<1, \mu_3>1$, then we have $\sigma_1=\sigma_2=-1$, and $\sigma_3=\sigma_4=1$, i.e., 
\[  f_1=(\mu_3^3  - \mu_4^3 )- \mu_2 ^3 (\mu_3^3 +  \mu_4^3 )=0, \
f_2=(\mu_2^3  +  \mu_4^3 )+\mu_3 ^3 (\mu_2^3  - \mu_4^3 )=0. \]
Then
\[ 0=f_{1}+f_{2}= \left(1-\mu_4^3\right) \left(\mu_2^3+\mu_3^3\right),   \]
thus $\mu_4=1$. By  the equation   $f_{1}=0$, we get $\mu_2^3=\frac{\mu_3^3-1}{\mu_3^3+1}$.  By  the third  equation, we have 
\[ 0=\frac{\left(\mu_3^3-1\right)^2 \mu_3^{18}}{\left(\mu_3^3+1\right)^2}+\left(\frac{\left(\mu_3^2-1\right) \left(\mu_3^3-1\right)^2}{\left(\mu_3^3+1\right)^2}-\mu_3^6\right)^3. \]
It is easy to verify that the above equation has no  root for   $\mu_3\in(1, \infty)$.
Thus, the system has no positive solutions. \\

\emph{Computations in the proof of Theorem \ref{thm:main-mnn011}}：

By assumption, the masses satisfy \eqref{IV123} \eqref{Tri124} \eqref{Tri134}, i.e.,  belong to the set 
\begin{align*}
\mathcal{V}_{IV}[123]  \bigcap&\left( \mathcal{V}_{III}[124] \bigcup  \mathcal{V}_{IV}[124]  \bigcup  \mathcal{V}_{V}[1,24]   \bigcup \mathcal{V}_{V}[2,14]\bigcup \mathcal{V}_{V}[4,12] \right) \\
 \bigcap&\left( \mathcal{V}_{III}[134] \bigcup  \mathcal{V}_{IV}[134]  \bigcup  \mathcal{V}_{V}[1,34]   \bigcup \mathcal{V}_{V}[3,14]\bigcup \mathcal{V}_{V}[4,13] \right). 
\end{align*}
We need to show that any  masses of the above set  is  not consistent with  any one of 
\[  \eqref{DiaIandII}, \  \eqref{III134}, \ \eqref{III234}, \ \eqref{IV234},  \]
and if the relation \eqref{IV134} is assumed, the masses are not consistent with any one of 
\[  \eqref{DiaIandIIagain}, \  \eqref{III124}.\]

We separate the  masses satisfying  \eqref{IV123} \eqref{Tri124} \eqref{Tri134} into the union of the following five subsets:
\begin{align*}
&{\it Subset \ 1}, \ &\mathcal{V}_{IV}[123]  \bigcap\left(  \mathcal{V}_{IV}[124]  \bigcup  \mathcal{V}_{V}[1,24]   \bigcup \mathcal{V}_{V}[2,14]\bigcup \mathcal{V}_{V}[4,12] \right) \\
& \  &\bigcap\left(  \mathcal{V}_{IV}[134]  \bigcup  \mathcal{V}_{V}[1,34]   \bigcup \mathcal{V}_{V}[3,14]\bigcup \mathcal{V}_{V}[4,13] \right)\\
 &{\it Subset  \ 2}, \ &(\mathcal{V}_{IV}[123]  \mathcal{V}_{III}[124] \mathcal{V}_{IV}[134] )  \bigcup (\mathcal{V}_{IV}[123]  \mathcal{V}_{IV}[124] \mathcal{V}_{III}[134] )  \\
  &{\it Subset  \ 3}, \ &\mathcal{V}_{IV}[123]  \mathcal{V}_{III}[124] \mathcal{V}_{III}[134] \\
   &{\it Subset  \ 4}, \ &\mathcal{V}_{IV}[123]  \mathcal{V}_{III}[124] \left( \mathcal{V}_{V}[1,34]   \bigcup \mathcal{V}_{V}[3,14]\bigcup \mathcal{V}_{V}[4,13]   \right)\\
  &{\it Subset  \ 5}, \ &\mathcal{V}_{IV}[123]  \mathcal{V}_{III}[134] \left( \mathcal{V}_{V}[1,24]   \bigcup \mathcal{V}_{V}[2,14]\bigcup \mathcal{V}_{V}[4,12]   \right). 
\end{align*}

For the first two subsets, straightforward computation show that there are only five positive solutions
\begin{align*}
&\left\{\mu_1=1,\mu_2=2,\mu_3=\sqrt{5},\mu_4=2\right\},   &\left\{\mu_1=1,\mu_2=\sqrt{5},\mu_3=2,\mu_4=2\right\},  \\
&\left\{\mu_1=1,\mu_2=2,\mu_3=\sqrt{3},\mu_4=2\right\},   &\left\{\mu_1=1,\mu_2=\sqrt{3},\mu_3=2,\mu_4=2\right\},  \\
&\left\{\mu_1=1,\mu_2=1/\sqrt{2},\mu_3=1/\sqrt{2}, \mu_4=1/\sqrt{2}\right\}. & 
\end{align*}
It is easy to check that none of the above  five solutions of  $\mu_i$'s  are  consistent with any one of $\eqref{DiaIandII},  \eqref{III134},  \eqref{III234},  \eqref{IV234}, $ 
and that those satisfy equation \eqref{IV134} do not satisfy  any one of 
$  \eqref{DiaIandIIagain},   \eqref{III124}.$

For the third subset, $\mathcal{V}_{IV}[123]  \mathcal{V}_{III}[124] \mathcal{V}_{III}[134]$,  the masses satisfy the system 
\[  	1+\sigma_1 \mu_2^2 +\sigma_2 \mu_3^2=0,\ 
1+\sigma_3/\mu_2 +\sigma_4/\mu_4=0, \ 
1+\sigma_5/\mu_3 +\sigma_6/\mu_4=0.  \]

 The above  system reduces  to  27 quartic equations of $\mu_4$, of which we can find their exact solutions.  There are  in total  eight positive  solutions. We will not write them explicitly here.  It is easy to check that the corresponding  masses are not consistent with \eqref{DiaIandII}, and  satisfy neither  $m_1+m_3+m_4=0,$ nor $m_2+m_3+m_4=0$. \\

The last two subsets are similar. We consider only the fourth one. 
Firstly, we claim that the fourth subset, $\mathcal{V}_{IV}[123]  \mathcal{V}_{III}[124] \left( \mathcal{V}_{V}[1,34]   \bigcup \mathcal{V}_{V}[3,14]\bigcup \mathcal{V}_{V}[4,13]   \right)$,  of masses are not consistent with  \eqref{DiaIandII}.

\textbf{Case I,  in $\mathcal{V}_{IA}[13,24] \bigcup \mathcal{V}_{IA}[14,23]$. }
If the masses are also in $\mathcal{V}_{IA}[13,24]$, then we have $m_1=m_3=1$, $m_2=m_4=-2$. Hence, the masses are not in $\mathcal{V}_{III}[124]$, a contradiction. Similarly, the masses can not belong to $\mathcal{V}_{IA}[14,23]$.

\textbf{Case II,   in $\mathcal{V}_{IB}[13,24]$ or in   $\mathcal{V}_{IB}[14,23]$. } 
We consider only the first one, i.e.,  the masses are also in $\mathcal{V}_{IB}[13,24]$.  Then $m_1m_3<0$, so the subset reduces to 
\[  \mathcal{V}_{IV}[123]  \mathcal{V}_{III}[124] \left( \mathcal{V}_{V}[1,34]   \bigcup \mathcal{V}_{V}[3,14]   \right)\mathcal{V}_{IB}[13,24].  \] 

\textbf{Subcase I, in $\mathcal{V}_{V}[1,34]$. } In this case, $m_3m_4>0$, so $m_1, m_2>0$, $m_3, m_4<0$ since $m_1+m_2+m_3=0$.  Then we have 
\begin{equation}\label{Subset4-3-1-1}
1+\mu_2^2-\mu_3^2=0. 
\end{equation}
So $1<\mu_3, \mu_2<\mu_3$.  The polynomial corresponding to 
$\mathcal{V}_{V}[1,34]$ is 
\begin{equation}\label{Subset4-3-1-2}
(\mu_3^2 +\mu_4^2)^2-\mu_3^3 \mu_4^3=0. 
\end{equation}
 Note that the relation corresponding to $\mathcal{V}_{III}[124]$ can be written as 
\[  (\mu_2-\frac{\mu_4}{1-\mu_4})(\mu_2-\frac{\mu_4}{1+\mu_4}) (\mu_2-\frac{\mu_4}{\mu_4-1}).      \]
There are three possibilities.

If $\mu_2=\frac{\mu_4}{1-\mu_4}$, then $\mu_4<\mu_2$.  Note that the polynomial \eqref{Subset4-3-1-1} becomes 
\begin{equation}\label{Subset4-3-1-3}
-\mu_3^2 \mu_4^2+2 \mu_3^2 \mu_4-\mu_3^2+2 \mu_4^2-2 \mu_4+1=0, 
\end{equation}
and the  polynomial corresponding to 
$\mathcal{V}_{IB}[13,24]$ is $-\mu_2^3(1-\mu_3^3)=\mu_4^3(1+\mu_3^3)$, or 
\begin{equation}\label{Subset4-3-1-4}
\left(\mu_3^3-1\right)-\left(\mu_3^3+1\right) (1-\mu_4)^3=0. 
\end{equation}
However, straightforward computation shows that \eqref{Subset4-3-1-2}, \eqref{Subset4-3-1-3}, \eqref{Subset4-3-1-4} has no common root.  Hence, the set is empty if $\mu_2=\frac{\mu_4}{1-\mu_4}$. 

Similarly, we can show that the set is empty if $\mu_2=\frac{\mu_4}{1+\mu_4}$,  or $\mu_2=\frac{\mu_4}{\mu_4-1}$.

\textbf{Subcase II, in $\mathcal{V}_{V}[3,14]$. } 
The set 
$\mathcal{V}_{IV}[123]  \mathcal{V}_{III}[124] \mathcal{V}_{V}[3,14]   \mathcal{V}_{IB}[13,24]$
is also empty. The proof is similar to the above case and is omitted. 

 To summarize,  the set $ \mathcal{V}_{IV}[123]  \mathcal{V}_{III}[124] \left( \mathcal{V}_{V}[1,34]   \bigcup \mathcal{V}_{V}[3,14]   \right)\mathcal{V}_{IB}[13,24]$ is empty.

\textbf{Case III,  in $\mathcal{V}_{II}[12,34]$. }
If the masses are also in $\mathcal{V}_{II}[12,34]$, then $\mu_2=\mu_3\mu_4.$  Then after the substitution $\mu_2=\mu_3\mu_4$,   the constraints corresponding to 
\[\mathcal{V}_{IV}[123], \  \mathcal{V}_{III}[124], \ \mathcal{V}_{V}[1,34]   \bigcup\mathcal{V}_{V}[3,14]   \bigcup  \mathcal{V}_{V}[4,13]   \]
 become three polynomials in terms of $(\mu_3,  \mu_4)$ with degree 12, 6 and 18. 
Straightforward computation shows that the  three polynomials have no common root, a contradiction. \\

Secondly, the fourth subset, $\mathcal{V}_{IV}[123]  \mathcal{V}_{III}[124] \left( \mathcal{V}_{V}[1,34]   \bigcup \mathcal{V}_{V}[3,14]\bigcup \mathcal{V}_{V}[4,13]   \right)$,  of masses satisfy neither $m_1+m_3+m_4=0,$ nor $m_2+m_3+m_4=0$. For instance, if  $m_1+m_3+m_4=0$, then $m_2=m_4$. The fact that the masses are in $\mathcal{V}_{III}[124] \mathcal{V}_{IV}[123]$ implies that the masses are 
  \[   \left\{ m_1=1, m_2=4, m_3=-5, m_4=4\right\}  \mbox{  or }    \left\{ m_1=1, m_2=-4, m_3=3, m_4=-4 \right\}.    \]
It is easy to check that both of them are not in the set $\mathcal{V}_{V}[1,34]   \bigcup \mathcal{V}_{V}[3,14]\bigcup \mathcal{V}_{V}[4,13]$.  \\

\emph{Computations in the proof of Theorem \ref{thm:main-mn012}}：

\textbf{Case 1}:  The masses are in  $\mathcal{V}_{I0}[13,24]\mathcal{V}[13,24] \mathcal{V}_{IV}[123]$. Recall that 
\[    V[13,24]=(\mathcal{V}_{I0}[12,34])\bigcup (\mathcal{V}_{I0}[14,23])\bigcup(\mathcal{V}_{II0}[13,24])\bigcup \mathcal{V}_{IV}[124]\bigcup \mathcal{V}_{IV}[234].  \]  
We ignore the first two choices $\mathcal{V}_{I0}[12,34]\bigcup \mathcal{V}_{I0}[14,23]$, since we have showed 
\[ \mathcal{V}_{I0}[13,24](\mathcal{V}_{I0}[12,34]\bigcup \mathcal{V}_{I0}[14,23]) =\emptyset,  \]
see Case 4 of the proof of  Theorem \ref{thm:main-63}. Hence, there are three possibilities. 

If the masses belong to $\mathcal{V}_{I0}[13,24]\mathcal{V}_{II0}[13,24] \mathcal{V}_{IV}[123]$. Then  $m_1m_3=m_2m_4$ and the masses are in $\mathcal{I}_{I}[13,24]$. Then we have  $m=0$, a contradiction.

If the masses belong to $\mathcal{V}_{I0}[13,24]\mathcal{V}_{IV}[124] \mathcal{V}_{IV}[123]$, then $m_3=m_4$. Set $m_3=s$, then $m_2=-(s+1)$. The fact that the masses are in $\mathcal{I}_{I}[13,24]$ implies 
\[  \frac{s}{\sqrt[3]{s+1} }  + s(s+1)=0, \]
which is impossible. Similarly, the masses can not belong to $\mathcal{V}_{I0}[13,24]\mathcal{V}_{IV}[234] \mathcal{V}_{IV}[123]$.

\textbf{Case 2}:  The masses are in  $\mathcal{V}_{II0}[12,34]\mathcal{V}[13,24] \mathcal{V}_{IV}[123]\mathcal{V}[23,14]$.  In fact, the set  $\mathcal{V}_{II0}[12,34]\mathcal{V}[13,24] \mathcal{V}_{IV}[123]$ is empty.  
Since  $m_2=m_3m_4$ and $1+m_2+m_3=0$, we have 
\[ m_1=1, \ m_2= \frac{-t}{1+t}, \  m_3= \frac{-1}{1+t}, \  m_4= t. \]
By Case 1, Diagram I can not occur, then
the masses are in 
$\mathcal{V}_{II0}[13,24]\bigcup \mathcal{V}_{IV}[124]\bigcup \mathcal{V}_{IV}[234].$

If the masses belong to $\mathcal{V}_{II0}[13,24]$, then $m_3=m_2m_4$, so we have $t=1$ (since $m_1+m_4\ne 0$). Thus, $m_2=m_3=-\frac{1}{2}$. This contradicts with the fact that the masses are in $\mathcal {I}_{II}[13,24]$, i.e., 
\[ 0= \frac{m_1+m_3}{\sqrt[3]{m_1+m_3}} +  \frac{m_3+m_2}{\sqrt[3]{m_3+m_2}} + \frac{m_2+m_4}{\sqrt[3]{m_2+m_4}} + \frac{m_4+m_1}{\sqrt[3]{m_4+m_1}}.    \]

If  the masses belong to $\mathcal{V}_{IV}[124]$, then $m_3=m_4$, or $(t+1)t+1=0$, which has no real solution. Similarly, the masses can not belong to $\mathcal{V}_{IV}[234]$. \\

\emph{Computations in the proof of Theorem \ref{thm:main-m01}}:

We want to show that the system 
\begin{align} \label{sys:thm81}
&\mu_2^3(  1+ \sigma_1  \mu_3^3 )= \sigma_2\mu_4^3(  1- \sigma_1  \mu_3^3 ) \\
&\mu_2^3(  1+ \sigma_3  \mu_4^3 )= \sigma_4\mu_3^3(  1- \sigma_3  \mu_4^3 ) \notag \\
&1+\mu_2^2 =\mu_3^2+\mu_4^2 \notag 
\end{align} 
has no solution, if we assume
$ (\mu_2-1)(\mu_3-1)(\mu_4-1)(\mu_2-\mu_3)(\mu_2-\mu_4)(\mu_3-\mu_4) \ne 0.$

Without lose of generality, assume $\mu_3>\mu_4$. There are three cases: $\mu_2>\mu_3>\mu_4$,  $\mu_3>\mu_4>\mu_2$,  and $\mu_3>\mu_2>\mu_4$. 

\textbf{Case 1}: $\mu_2>\mu_3>\mu_4$.  Then  $\mu_4>1, \ \sigma_1=\sigma_3=-1, \  \sigma_2=\sigma_4=-1. $
Therefore, we have 
\[  \mu_3^3 \frac{1- \mu_3^3}{1+\mu_3^3}=\mu_4^3 \frac{1-  \mu_4^3}{1+\mu_4^3},  \]
which has no solution since $t^3 \frac{1- t^3}{1+t^3}$ is an increasing function on $(1, \infty)$. 

\textbf{Case 2}: $\mu_3>\mu_4>\mu_2$. System \eqref{sys:thm81} has no solution. The proof is similar to the above case and is omitted.

\textbf{Case 3}: $\mu_3>\mu_2>\mu_4$.  Note that $\mu_4<1<\mu_3$. Rewrite the first two equation as 
\begin{align*}
&(\mu_2^3+ \sigma_1  \mu_4^3 )= \sigma_2\mu_3^3(\mu_2^3- \sigma_1  \mu_4^3 ),\ \mu_2^3(    1+ \sigma_3  \mu_4^3 )= \sigma_4\mu_3^3(  1- \sigma_3  \mu_4^3 ). 
\end{align*} 
Then it is easy to see that $\sigma_1=\sigma_2=\sigma_3=\sigma_4=1.$ Therefore, we have 
\begin{equation} \label{equ:thm81}
\frac{\mu_2^3+ \mu_4^3}{\mu_2^3(\mu_2^3- \mu_4^3)}=\frac{1+ (\mu_4/\mu_2)^3}{\mu_2^3(1- (\mu_4/\mu_2)^3)}=\frac{1+ \mu_4^3}{1-\mu_4^3}. 
\end{equation}
Note that $\mu_4, \mu_4/\mu_2 \in (0,1)$ and $\frac{1+ t^3}{1-t^3}$ is an increasing function on $(0,1)$.  If $\mu_2<1$, then $\mu_4<\mu_4/\mu_2<1$, and 
\[  \frac{1+ (\mu_4/\mu_2)^3}{\mu_2^3(1- (\mu_4/\mu_2)^3)} >\frac{1+ (\mu_4/\mu_2)^3}{(1- (\mu_4/\mu_2)^3)} >\frac{1+ \mu_4^3}{1-\mu_4^3},   \]
a contradiction.  Similarly, equation \eqref{equ:thm81} can not hold if   $\mu_2>1$.

\end{document}